%% file: zgammacircunique.tex
\documentclass[a4paper,10pt]{article}

\pdfoutput=1

\usepackage{amsmath}        % use AMS packages
\usepackage{amssymb}
\usepackage{amsthm}
\usepackage[nooneline,hang]{subfigure}
\usepackage[font=normalsize,singlelinecheck=off,labelfont=bf]{caption}
\usepackage{graphicx}
\usepackage{bbold}
\usepackage[usenames]{color}
\usepackage{dsfont}
\usepackage{floatflt}
\usepackage{enumerate}
\usepackage{mathrsfs}
\usepackage{url}

\usepackage{hyperref}

\usepackage[a4paper,hmargin={1.5in,1.5in},vmargin={1.5in,1.5in}]{geometry}
\usepackage{setspace}

\newcommand{\R}{\mathds R}
\newcommand{\C}{\mathds C}
\newcommand{\N}{\mathds N}
\newcommand{\Z}{\mathds Z}
\newcommand{\Sp}{\mathds S}
\newcommand{\V}{\mathds V}

\newcommand{\eps}{\varepsilon}
\renewcommand{\Re}{\text{Re}}
\renewcommand{\Im}{\text{Im}}

\newtheoremstyle{plainNoItalics}{}{}{\normalfont}{}{\bfseries}{.}{ }{}

\theoremstyle{plain}
\newtheorem{theorem}{Theorem}[section]
\newtheorem{proposition}[theorem]{Proposition}
\newtheorem{corollary}[theorem]{Corollary}
\newtheorem{lemma}[theorem]{Lemma}

\theoremstyle{plainNoItalics}
\newtheorem{remark}[theorem]{Remark}
\newtheorem{example}[theorem]{Example}

\newtheorem*{theorem*}{Theorem}
\newtheorem*{proposition*}{Proposition}
\newtheorem*{lemma*}{Lemma}
\newtheorem*{corollary*}{Corollary}
\newtheorem*{remark*}{Remark}

\newtheorem*{observation*}{Observation}
\newtheorem*{example*}{Example}
\newtheorem*{examples*}{Examples}
\newtheorem*{assumption*}{Assumption}

\newtheorem*{acknowledgement}{Acknowledgment}

\theoremstyle{definition}
\newtheorem{definition}[theorem]{Definition}
\newtheorem*{definition*}{'Definition'}

\newtheorem*{definitionu*}{Definition}

\title{Rigidity of quasicrystallic and $Z^\gamma$-circle patterns}
\author{Ulrike B\"ucking\thanks{Institut f\"ur Mathematik, Technische
Universit\"at Berlin, Stra\ss{}e des 17.~Juni 136, 10623 Berlin, Germany.
E-mail: {\tt buecking@math.tu-berlin.de}.
}
}

\begin{document}

%%% Achtung !!!!!!!!! Fuer bessere Bilder ``3'' am Ende der Dateinamen
%%% entfernen !!!!!!!!!!!!!!!!!

\maketitle

\begin{abstract}
The uniqueness of the orthogonal $Z^\gamma$-circle patterns as studied by
Bobenko and Agafonov is shown, given the combinatorics and some boundary
conditions.
Furthermore we study (infinite) rhombic embeddings in the plane which are
quasicrystallic, that is they have only finitely many
different edge directions. Bicoloring the vertices of the
rhombi and adding circles with centers at vertices of one of the colors and
radius equal to the edge length leads to isoradial quasicrystallic circle
patterns.
We prove for a large class of such circle patterns which
cover the whole plane that they are uniquely determined up to affine
transformations by the combinatorics and the intersection angles. Combining
these two results, we obtain the rigidity of large classes of quasicrystallic
$Z^\gamma$-circle patterns.
\end{abstract}

\section{Introduction}

Circles, especially circle packings and circle patterns, have
successfully been used over the past years to define and study discrete analogs
of classical smooth objects. In particular, this approach leads to discrete
holomorphic mappings, for example discrete analogs of the power functions
$z^\gamma$, and to discrete holomorphic function theory.
See for example \cite{Sch97,Bo99,BH03,BS08} for some of the contributions to
the theory of circle patterns and
\cite{St05} for results on circle packings.

In this article we focus on circle patterns which are
characterized by a given combinatorics specifying which circles should intersect
and by the corresponding intersection angles. Thus we associate to a circle
pattern a pattern of kites corresponding to intersecting circles, see
Figures~\ref{figbquad} (right) and~\ref{figPenrose}.
A particularly suitable source for the required knowlegde on circle patterns
and their relations to consistency, integrability, and discrete
holomorphic functions is the textbook~\cite{BS08} in discrete differential
geometry.
Our main results can roughly be summarized as follows:
The combinatorics (together with suitable intersection angles and boundary
conditions, if necessary) determines the geometry of the circle pattern.
This rigidity result can be interpreted as a discrete version of Liouville's
Theorem in complex analysis.

\begin{figure}[t]
\begin{center}
\includegraphics[height=4.cm]{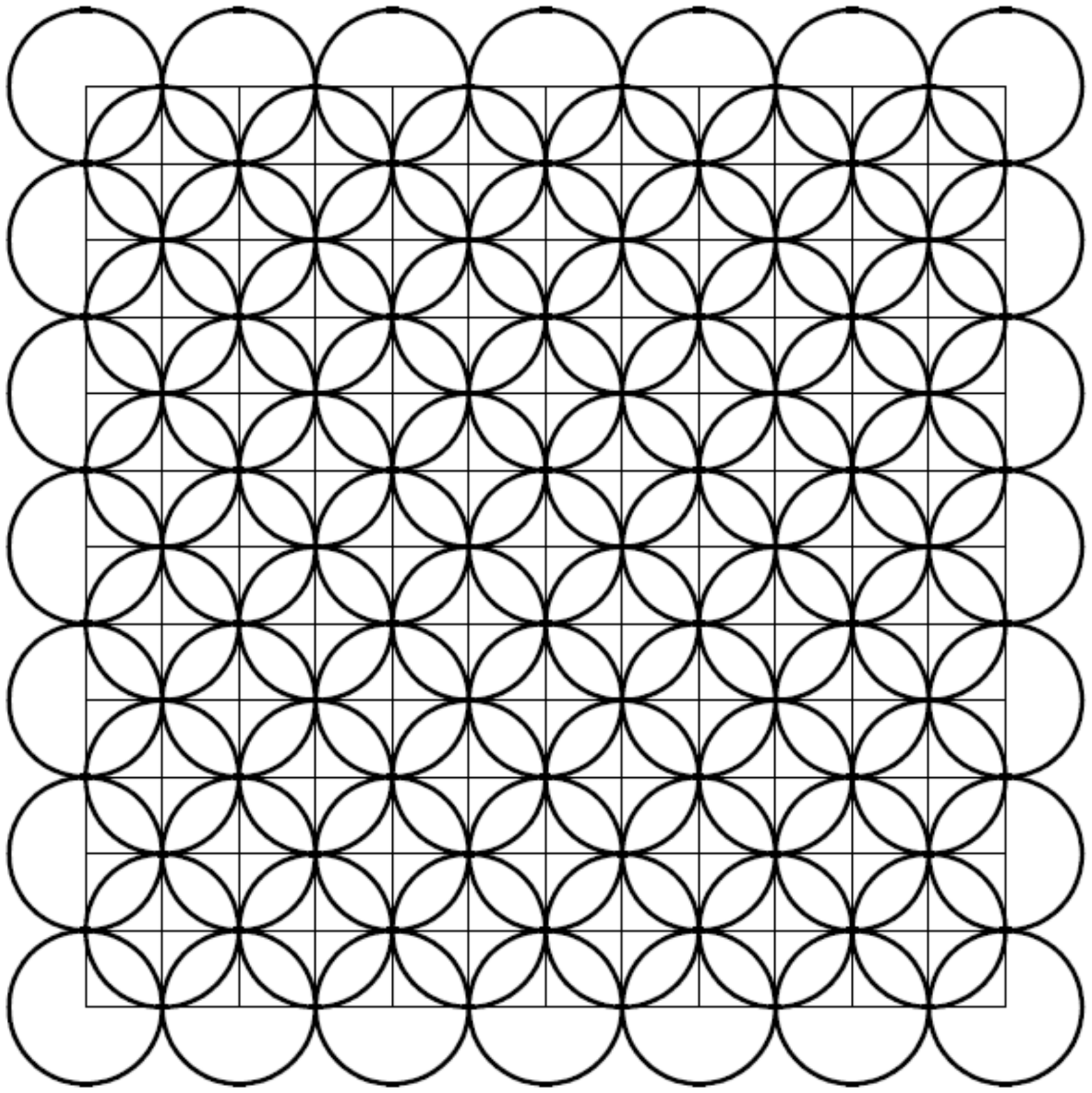}
\hspace{1cm}
\begin{picture}(25,75)
 \put(5,37){$\longrightarrow$}
\end{picture}
\hspace{1cm}
\includegraphics[height=3cm]{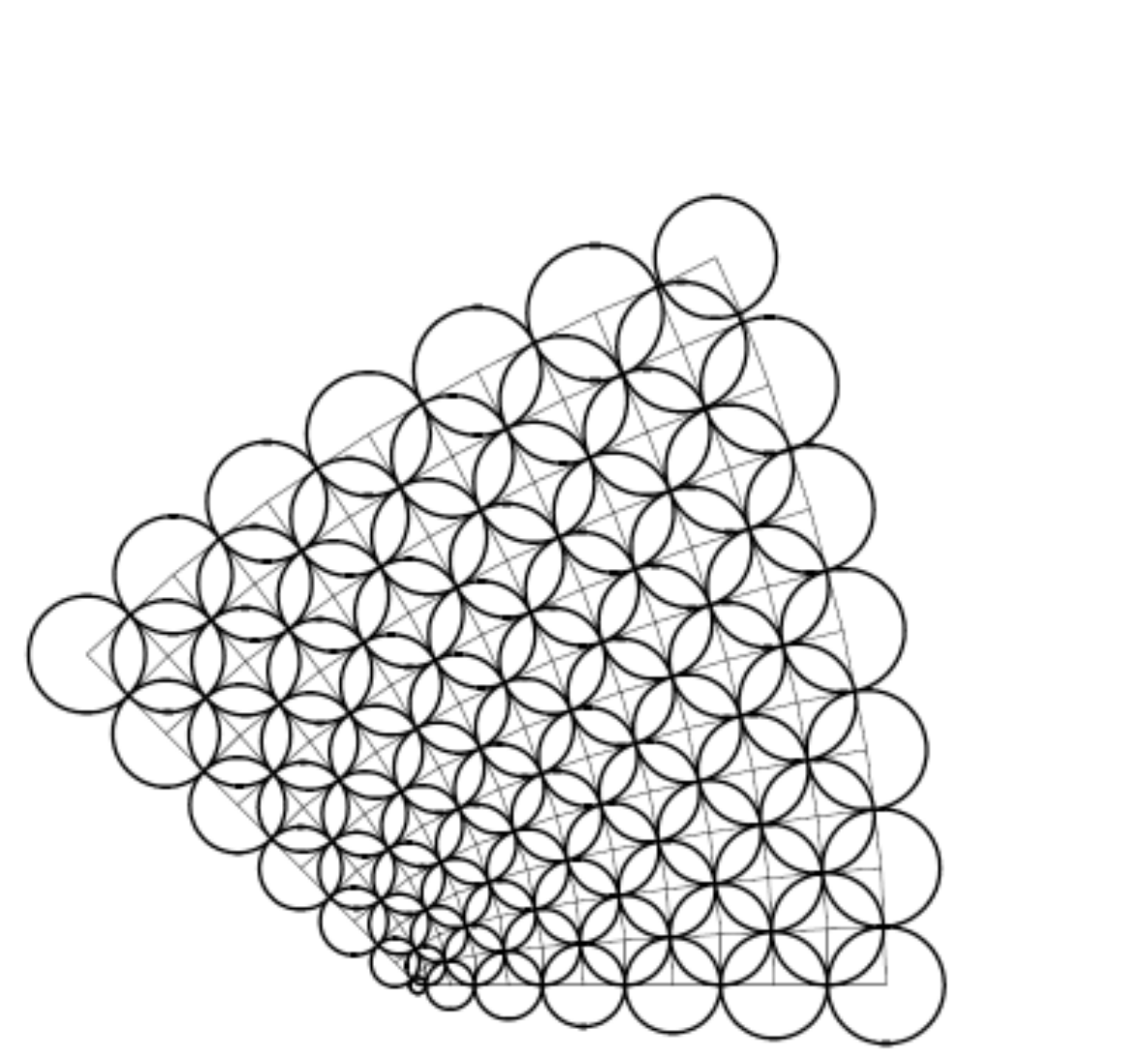}
\end{center}
\caption{Illustration of the discrete conformal map $Z^{3/2}$ and the orthogonal
$Z^{3/2}$-circle pattern ({\it right}).}\label{figZgamma}
\end{figure}

Our first result concerns the uniqueness of
circle patterns which are discrete analogs of the power functions $z^\gamma$
for $\gamma\in(0,2)$ as defined in~\cite{Bo99,AB00,Ag05}.
Here, the square grid combinatorics, the orthogonal intersection angles, and
the two boundary lines of a sector uniquely determine the geometry of the
pattern, see Figure~\ref{figZgamma} for an illustration.
The rigidity of orthogonal $Z^\gamma$-circle patterns
was only known for rational $\gamma$, see~\cite{Ag05}.

Furthermore, we consider the case of a circle pattern which covers the whole
 complex plane and for which the radii of all circles are equal and
the interiors of different kites are disjoint. Then the corresponding kites
form a rhombic embedding. Also, assume that there are only finitely many
different edge directions of the kites. Such rhombic embeddings are called
{\em quasicrystallic}~\cite{BMS05}.
Additionally, orient an edge $\vec{e}$ and consider a line perpendicular to
the edge. Move this line parallelly in positive and in negative direction along
$\vec{e}$. We suppose that in both cases this moving line intersects infinitely
often edges parallel to $\vec{e}$. We assume that this property is true at least
for two edges with linearly independent directions. Such rhombic embeddings of
the plane are for example given by Penrose tilings, see for example
Figure~\ref{figPenrose}. We show that
any other embedded circle pattern with the same combinatorics and intersection
angles of the kites is the image of this embedding by an affine transformation.

Rigidity of some classes of infinite circle patterns of the plane have already
been studied. Schramm~\cite{Sch97} considers square grid combinatorics and
orthogonal intersection angles. As an essential step of our proof we generalize
his result to square grid circle patterns with regular intersection angles
$\psi\in(0,\pi)$ and $(\pi-\psi)$.
He~\cite{He99} studies disk triangulation graphs and exterior intersection
angles in~$[\pi/2,\pi]$ which does not cover the class of isoradial
quasicrystallic circle patterns defined above.

As observed in~\cite{BMS05} isoradial quasicrystallic rhombic embeddings can be
used to define corresponding quasicrystallic $Z^\gamma$-circle patterns, see
also~\cite{BS08}.
Examples of quasicrystallic $Z^\gamma$-circle patterns are shown
in Figures~\ref{figVeronika}, \ref{figVeronika2}, and~\ref{figVeronika3}. They
have been created using software developed by Veronika Schreiber for her
diploma thesis~\cite{Veronika}.
Our rigidity result for orthogonal $Z^\gamma$-circle patterns is also
generalized for large classes of quasicrystallic $Z^\gamma$-circle patterns.

There is some more literature concerning
rigidity for infinite planar circle {\em packings}, that is configurations with
non-overlapping touching circles, see~\cite{RS87,Du97,Sch91,HSch93}.
Our rigidity proofs adapt some of the ideas which have been used for
packings. In particular, we apply discrete potential theory.

This paper is organized as follows.
First we introduce terminology and present useful facts about circle patterns.
We especially focus on regular circle patterns
with square grid combinatorics. Then we recall in Section~\ref{secZgammaOrth}
the definition and
some properties of the orthogonal $Z^\gamma$-circle patterns for $\gamma\in (0
,2)$ as studied in~\cite{AB00,Ag03,Ag05} and prove their rigidity. In
Section~\ref{secQuasiCirc} we introduce quasicrystallic circle patterns and
prove uniqueness for a class of these patterns. Finally, we recall the
definition and some facts about quasicrystallic $Z^\gamma$-circle patterns for
$\gamma\in (0 ,2)$ as studied in~\cite{Agpsi,BMS05,BS08} and prove rigidity
for certain classes of these patterns. A more detailed version of the results
can be found in~\cite{diss}.

%%%%%%%%%%%%%%%%%%%%%%%%%%%%%%%%%%%%%%%%%%%%%%%%%%%%%%%%%%%%%%%%%%%%%%%%%%%
\section{Circle patterns}\label{CirclePatterns}
In this section we focus on a definition and some useful properties of circle
patterns. We describe circle patterns using combinatorial data and
intersection angles.

The combinatorics are specified by a {\em b-quad-graph} $\mathscr D$, that
is a strongly regular cell
decomposition of a domain in $\C$ possibly with boundary such that all
2-cells (faces) are embedded  and counterclockwise oriented.
Furthermore all faces of $\mathscr D$ are
quadrilaterals, that is there are exactly four edges incident to
each face, and the 1-skeleton of $\mathscr D$ is a
bipartite graph.
We always assume that the
vertices of $\mathscr D$ are colored white and black. To these two sets of
vertices we associate two planar graphs $G$ and $G^*$ as follows.
The vertices $V(G)$ are all white vertices of 
$V(\mathscr D)$. The edges $E(G)$ correspond to faces of $\mathscr D$, that is
two vertices of $G$ are connected by an edge if and only if they
are incident to the same face. The
dual graph $G^*$ is constructed analogously by taking as vertices
$V(G^*)$ all black vertices of $\mathscr D$.
$\mathscr D$ is called {\em simply
  connected} if it is the cell decomposition of a simply connected
domain of $\C$
and if every closed chain of faces is null homotopic in $\mathscr D$.

For the intersection angles, we use a labelling
$\alpha:F({\mathscr D})\to (0,\pi)$ of the faces of $\mathscr D$. By abuse of
notation, $\alpha$ can also be understood as a function defined on $E(G)$ or on
$E(G^*)$.
The labelling $\alpha$ is called {\em admissible} if it satisfies the following
condition at all interior black vertices $v\in V_{int}(G^*)$:
\begin{equation}\label{condalpha}
  \sum_{f \text{ incident to } v} \alpha(f) = 2\pi.
\end{equation}

\begin{figure}[tb]
\begin{center}
 \includegraphics[width=0.45\textwidth]{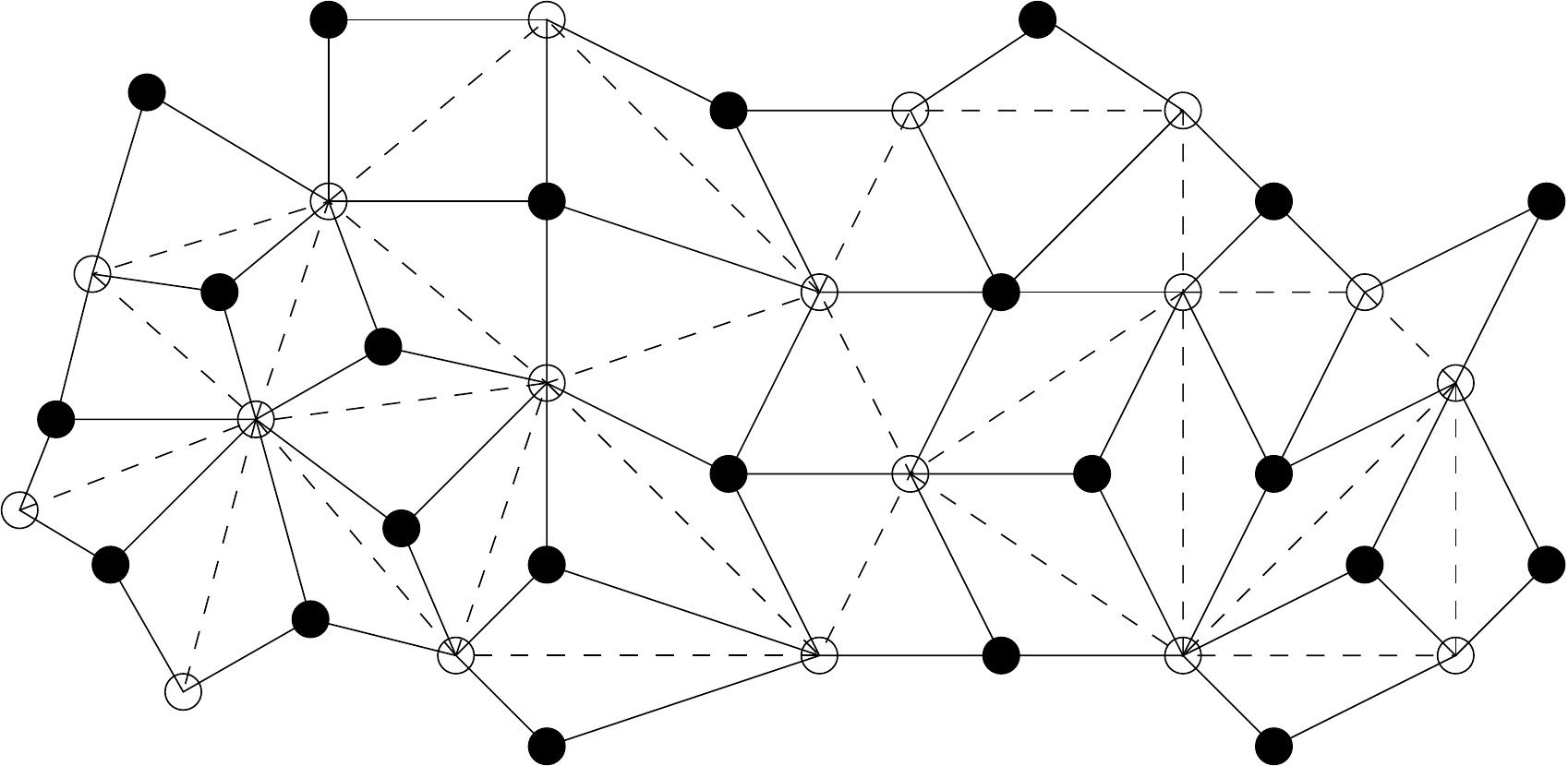}
\hspace{2em}
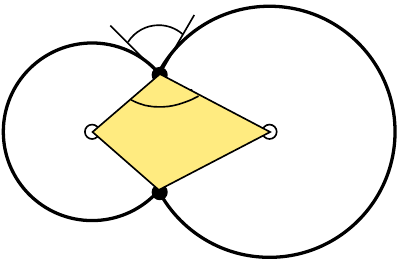
\end{center}
\caption{{\it Left:} An example of a b-quad-graph $\mathscr D$ (black edges and
bicolored
  vertices) and its associated graph $G$ (dashed edges and white vertices).
{\it Right:} The exterior intersection angle $\alpha$
of two intersecting circles and the associated
  kite built from centers and intersection points.
$\beta=\pi-\alpha$ is the interior intersection angle.}\label{figbquad}
\end{figure}

\begin{definition}\label{defcircpattern}
Let $\mathscr D$ be a b-quad-graph with associated graph $G$ and let
$\alpha:E(G)\to (0,\pi)$
be an admissible labelling.
  An {\em (immersed planar) circle pattern for  $\mathscr D$ (or $G$) and
    $\alpha$} are an indexed collection ${\mathscr 
    C}=\{C_z:z\in V(G)\}$ of circles in $\C$ and an indexed
  collection ${\mathscr K}=\{K_e:e\in E(G)\}=\{K_f:f\in F({\mathscr
    D})\}$ of closed
  kites, which all carry the same orientation, such that the following
conditions hold.
\begin{enumerate}[(1)]
  \item\label{intersectionPoint}
If $z_1,z_2\in V(G)$ are incident vertices in $G$,
the corresponding circles $C_{z_1},C_{z_2}$ intersect with exterior
intersection angle $\alpha([z_1,z_2])$.
 Furthermore, the kite $K_{[z_1,z_2]}$
    is bounded by the centers of the circles $C_{z_1},C_{z_2}$,
    the two intersection points, and the corresponding edges, as in
Figure~\ref{figbquad} (right).
 The intersection points
are associated to black vertices of $V(\mathscr D)$ or to vertices of $V(G^*)$.
\item If two faces are incident in
    $\mathscr D$, the corresponding kites share a common edge.
  \item Let $f_1,\dots,f_n\in F({\mathscr D})$ be the faces
    incident to an interior vertex $v\in V_{int}({\mathscr D})$. Then 
    the kites $K_{f_1},\dots,K_{f_n}$ have mutually disjoint
    interiors.
The union $K_{f_1}\cup\dots\cup
    K_{f_n}$ is homeomorphic to a closed disk and contains the point $p(v)$
corresponding to $v$ in its interior.
  \end{enumerate}
  The circle pattern is called {\em embedded} if all kites of
  $\mathscr K$ have mutually disjoint interiors.
 The circle pattern is called {\em isoradial} if all
  circles of $\mathscr C$  have the same radius.
\end{definition}

  Note that we associate a circle pattern $\mathscr C$ to an immersion of the
kite pattern $\mathscr K$ corresponding to ${\mathscr D}$ where the edges
incident to the same white vertex are of equal length. The kites may also be
  non-convex and
can be constructed from a (suitable) given set of circles and from the
combinatorics of $G$.

As all circle patterns considered in this article will be planar and immersed,
the notion ``circle pattern'' will include these properties in the following.

%%%%%%%%%%%%%%%%%%%%%%%%%%%%%%%%%%%%%%%%%%%%%%%%%%%%%%%%%%%
% \subsection{The radius function}

For our study of a circle pattern ${\mathscr C}$ we will use the
radius function $r_{\mathscr C}=r$ which assigns to every vertex $z\in V(G)$
the radius $r_{\mathscr C}(z)=r(z)$ of the corresponding circle $C_z$. The
index $\mathscr C$ will be dropped whenever there is no confusion likely.
The following proposition specifies a condition for a radius function
to originate from a planar circle pattern, see~\cite{BS02} for a proof.

\begin{proposition}\label{propRadius}
Let $G$ be associated to a b-quad-graph $\mathscr
      D$ and let $\alpha$ be an admissible labelling. 

      Suppose that ${\mathscr C}$ is a planar circle
      pattern for ${\mathscr D}$ and $\alpha$ with radius function
      $r=r_{\mathscr C}$. Then for 
      every interior vertex $z_0\in V_{int}(G)$ we have
      \begin{equation} \label{eqFgen}
        %F(r(z_0),r(z_0+1),r(z_0+i),r(z_0-1),r(z_0-i)) := \\
        \Biggl(\sum_{[z,z_0]\in E(G)} 
        f_{\alpha([z,z_0])}(\log r(z)-\log r(z_0))\Biggr) -\pi=0,
        \end{equation}
       where
        \[ f_\theta(x):=\frac{1}{2i}\log
        \frac{1-e^{x-i\theta}}{1-e^{x+i\theta}},\]
        and the branch of the logarithm is chosen such that
        $0<f_\theta(x)<\pi$.

Conversely, suppose that $\mathscr D$ is simply connected and
      that $r:V(G)\to(0,\infty)$ satisfies~\eqref{eqFgen} for
      every $z\in  V_{int}(G)$. Then there is a planar circle pattern
      for $G$ and $\alpha$
with radius function $r$.
This pattern is unique up to isometries of $\C$.
\end{proposition}

For the special case of orthogonal
circle patterns with the combinatorics of the square grid, there are also other
characterizations, see for example~\cite{Sch97}.

Note that $2f_{\alpha([z,z_0])}(\log r(z)-\log r(z_0))$ is the angle at $z_0$
of the kite
with edge lengths $r(z)$ and $r(z_0)$ and angle $\alpha([z,z_0])$, as in
Figure~\ref{figbquad} (right).
Equation~\eqref{eqFgen} is the closing
condition for the closed chain of kites which correspond to the edges incident
to $z_0$. This corresponds to condition (3) of
Definition~\ref{defcircpattern}.

For further use we mention some properties of $f_\theta$,
see~\cite[Lemma~2.2]{Spr03}.

\begin{lemma}\label{lemPropf}
  \begin{enumerate}[(i)]
    \item 
      The derivative of $f_\theta$ is
    $ f_\theta'(x)=\frac{\sin\theta}{2(\cosh x -\cos \theta)}>0$.
    \item The function $f_\theta$ satisfies the functional equation
      $f_\theta(x)+f_\theta(-x)=\pi-\theta$.
\end{enumerate}
\end{lemma}

If there exists an isoradial circle pattern, we can obtain another circle
pattern from a given radius function.

\begin{lemma}\label{lemInvRad}
 Let $G$ be a graph constructed from a b-quad-graph $\mathscr
      D$ and let $\alpha$ be an admissible labelling. 
      Suppose that there exists an isoradial circle pattern for $G$ and
$\alpha$. Let $r$ be the radius function of a planar circle
      pattern for ${\mathscr D}$ and $\alpha$. Then there is a circle
pattern $\mathscr C$ for $G$ and $\alpha$ with radius function $r_{\mathscr
C}=1/r$.
\end{lemma}
\begin{proof}
 By Lemma~\ref{lemPropf}~(ii) the function
$1/r$ satisfies condition~\eqref{eqFgen} of Proposition~\ref{propRadius} for
all interior vertices $z_0\in V_{int}(G)$. In particular,
\begin{multline*}
\sum_{[z,z_0]\in E(G)}
        f_{\alpha([z,z_0])}(\log \frac{1}{r(z)}-\log
\frac{1}{r(z_0)})\\
=\sum_{[z,z_0]\in E(G)} (\pi -\alpha([z,z_0])) - \sum_{[z,z_0]\in E(G)}
        f_{\alpha([z,z_0])}(\log \frac{r(z)}{r(z_0)}) 
=2\pi-\pi =\pi.
\end{multline*}
Here we have also used that $\sum_{[z,z_0]\in E(G)} (\pi
-\alpha([z,z_0]))
=2\pi$
since there is an isoradial circle pattern for $G$ and $\alpha$ and the
assumption that $r$ is the radius function of a circle pattern for $G$ and
$\alpha$.
\end{proof}

Let $G$ be a graph constructed from a b-quad-graph $\mathscr
      D$ and let $\alpha$ be an admissible labelling. %of the edges $E(G)$.
      Suppose that ${\mathscr C}_1$ and ${\mathscr C}_2$ are planar circle
      patterns for ${\mathscr D}$ and $\alpha$ with radius functions
      $r_1=r_{{\mathscr C}_1}$ and $r_2=r_{{\mathscr C}_2}$ respectively.
Define a {\em comparison function} $w:V({\mathscr D})\to \C$ by
\begin{equation}\label{eqdefw}
\begin{cases} w(y)=r_2(y)/r_1(y) &\text{for } y\in V(G),\\
w(x)=\text{e}^{i\delta(x)}\in \Sp^1 & \text{for } x\in V(G^*)
. \end{cases} 
\end{equation}
Here $\delta(x)\in\R$ or $w(x)=\text{e}^{i\delta(x)}$ is defined to be
the rotation angle or the rotation respectively of the edge-star at $x\in
V(G^*)$ when changing from
the circle pattern ${\mathscr C}_1$ to ${\mathscr C}_2$.
Note that $w(y)$ is the scaling factor of the circle corresponding to $y\in
V(G)$.
Then $w$ satisfies the following {\em Hirota Equation} for
all faces $f\in F({\mathscr D})$.
\begin{equation}\label{eqw}
  w(x_0)w(y_0)a_0 -w(x_1)w(y_0)a_1
  -w(x_1)w(y_1)a_0 +w(x_0)w(y_1)a_1 =0
\end{equation}
Here $x_0,x_1\in V(G^*)$ and $y_0,y_1\in V(G)$ are the black and white
vertices incident to $f$ and $a_0=x_0-y_0$ and $a_1=x_1-y_0$ are the
directed edges. Thus equation~\eqref{eqw} is the
closing condition for the kite corresponding to the face $f$.
Furthermore, the Hirota Equation is
3D-consistent; see Sections~10 and~11 of~\cite{BMS05} or~\cite{BS08} for
more details. This property will be used in Section~\ref{secQuasiCirc}.

\section{SG-circle patterns}

In this paper we are particularly interested in the special case of regular
circle patterns with square grid combinatorics.
First, we fix some notation.
Let $SGD$ be the regular square grid cell decomposition of the complex
plane, that is the vertices are $V(SGD)=\Z+i\Z$ and the edges are given by
pairs of vertices $[z,z']$ with $z,z'\in V(SGD)$ and $|z-z'|=1$.
The 2-cells are squares $\{z+a+ib: a,b\in[0,1]\}$ for $z\in
V(SGD)$. As $SGD$ is a b-quad-graph, the vertices
\[V(SG)=\{n+im\in\Z+i\Z :n+m=0\pmod{2}\} \]
are assumed to be colored white.
As above, $SG$ and its dual $SG^*$ are defined as the associated graphs to
$SGD$.
Furthermore,
$SG(n,v)$ with $n\in\N$ and $v\in V(SG)$ denotes the subgraph of all vertices
with combinatorial distance at most $n$ from $v$ in $SG$.

Let $\psi\in (0,\pi)$ be a fixed angle. Define the
following regular labelling $\alpha_\psi$ on $E(SG)$.
Let $[z_1,z_2]$ be an edge connecting the
vertices $z_1,z_2\in V(SG)$. Without loss of generality, we assume that
$\Re(z_1)< \Re (z_2)$. Then
\begin{equation}\label{defalphapsi}
\alpha_\psi([z_1,z_2])=\begin{cases}\psi &\text{if } \Im(z_1)<
\Im(z_2),\\
\pi-\psi &\text{if } \Im(z_1)> \Im(z_2).\end{cases}
\end{equation}

If $G$ is a subgraph of $SG$, a circle pattern for $G$ and
$\alpha_\psi$ is called {\em $SG$-circle pattern}.
The choice $\psi=\pi/2$ leads to orthogonal $SG$-circle patterns as
considered by Schramm in~\cite{Sch97}.

\begin{theorem}[Rigidity of $SG$-circle patterns]\label{theoSGUniq}
 Suppose that $\mathscr C$ is an embedded planar circle pattern for $SG$ and
$\alpha_\psi$. Then $\mathscr C$ is the image of a regular isoradial circle
pattern for $SG$ and $\alpha_\psi$ under a similarity.
\end{theorem}
The proof is a suitable adaption of the corresponding proof for orthogonal
$SG$-circle
patterns given by Schramm using suitable M\"obius invariants,
see~\cite[Theorem~7.1]{Sch97} or~\cite{diss}. This adaption needs
the following generalization of the Ring Lemma
of~\cite{RS87}, which is also useful in the following.

\begin{lemma}\label{lemRinggen}
 Let $r$ be the radius function of an embedded circle pattern for $SG(3,0)$ and
$\alpha_\psi$. There is a constant $C=C(\psi)>0$, independent of $r$, such
that for $k=0,1,2,3$ there holds
\begin{equation*}
 \frac{r(i^k(1+i))}{r(0)} > C .
\end{equation*}
\end{lemma}
\begin{proof}
 Assume the contrary. Then there is a sequence of embedded circle patterns for
$SG(3,0)$ and $\alpha_\psi$ such that $r_n(0)=1$ and $r_n(i^k(1+i))\to 0$ as 
$n\to\infty$ for
some $k\in\{0,1,2,3\}$. Without loss of generality we assume that $k=0$. We also
may assume that the circle $C_0$ corresponding to the vertex $0\in V(SG)$ and
the intersection point corresponding to $1\in V(SG^*)$ are fixed for the whole
sequence. Then there is a subsequence such that all the circles converge to
circles or lines, that is converge in the Riemann sphere $\hat{\C}\cong \Sp^2$.
Now equation~\eqref{eqFgen} implies that there exist some kites
which intersect in the limit in their interiors.
But this is a contradition to the embeddedness of the sequence.
\end{proof}

If the number of surrounding generations is big enough, there is the following
useful estimation on the quotient of radii of incident vertices.

\begin{lemma}\label{corSn}
There is an absolute constant $C>0$ such that the following holds.

Let $G$ be a subgraph of $SG$ and let $\mathscr C$ be
  an embedded circle pattern for $G$ and $\alpha_\psi$
with radius function $r$. Let $v\in V(G)$ be a vertex and suppose that
 $SG(n,v)\subset G$, that is
$G$ contains $n$ generations of $SG$ around
  $v$, for some $n\geq 3$. Then for all vertices $w$ incident to $v$
  there holds
  \begin{equation}
    1-\frac{C}{n} \leq \frac{r(w)}{r(v)}\leq 1+\frac{C}{n}.
  \end{equation}
\end{lemma}
The proof is an adaption of the corresponding proof for hexagonal
circle packings given by Aharonov in~\cite{Ah90,Ah94} and uses
Lemma~\ref{lemRinggen}, see~\cite{diss} for more details. The necessary results
on discrete potential theory can be found in the appendix
of~\cite{paperconv} and in~\cite{diss}.

%%%%%%%%%%%%%%%%%%%%%%%%%%%%%%%%%%%%%%%%%%%%%%%%%%%%%%%%%%%%%%%%%%%%
\section{Uniqueness of orthogonal $Z^\gamma$-circle
patterns}\label{secZgammaOrth}

An orthogonal
circle pattern with the combinatorics of the square grid associated to the map
$z^\gamma$ was introduced by Bobenko in~\cite{Bo99}. Further development of the
theory is due to Agafonov and Bobenko~\cite{AB00,Ag03,Ag05}.

\subsection{Definition and useful properties}
In the following we briefly summarize the definition and some known facts about
orthogonal $Z^\gamma$-circle patterns, see~\cite{AB00,Ag03,Ag05} for more
details.

\begin{definition}
  Let $D\subset \Z^2$. A map $f:D\to \C$ is called {\em discrete
    conformal} if all its elementary quadrilaterals are conformal
  squares, i.e.\ their cross-ratios are equal to $-1$:
\begin{multline}\label{defZgamma1}
q(f_{n,m}, f_{n+1,m}, f_{n+1,m+1}, f_{n,m+1}):=  \\
\frac{(f_{n,m}- f_{n+1,m})(f_{n+1,m+1}- f_{n,m+1})}{(f_{n+1,m}-
  f_{n+1,m+1})(f_{n,m+1}- f_{n,m})} = -1.
\end{multline}
Here and below we abbreviate $f_{n,m}=f(n,m)$.

A discrete conformal map $f_{n,m}$ is called {\em embedded} if the
interiors of different elementary quadrilaterals $(f_{n,m}, f_{n+1,m},
f_{n+1,m+1}, f_{n,m+1})$ are disjoint. 
\end{definition}

Note that the definition of a discrete conformal map is M\"obius invariant and
is motivated by
the following characterization for smooth mappings:\\
A smooth map $f:\C\subset D\to \C$ is called {\em conformal} (holomorphic or
antiholomorphic) if and only if for all $z=x+iy\in D$ there holds
\[ \lim_{\eps\to 0} q(f(x,y), f(x+\eps,y), f(x+\eps,y+\eps),
f(x,y+\eps))= -1.\]

In order to construct an embedded discrete analog of $z^\gamma$ the following
approach is used. Equation~\eqref{defZgamma1} can be
supplemented with the nonautonomous constraint
\begin{multline}\label{defZgamma2}
  \gamma f_{n,m}= 2n \frac{(f_{n+1,m}- f_{n,m})(f_{n,m}-
    f_{n-1,m})}{(f_{n+1,m}- f_{n-1,m})} \\
+ 2m \frac{(f_{n,m+1}-
    f_{n,m})(f_{n,m}- f_{n,m-1})}{(f_{n,m+1}- f_{n,m-1})}.
\end{multline} 
This constraint, as well as its compatibility with~\eqref{defZgamma1},
is derived from some monodromy problem; see~\cite{AB00}.
We assume that $0<\gamma<2$ and denote
\[\Z_+^2= \{(n,m)\in \Z^2: n,m\geq 0\}. \]
The asymptotics of the constraint~\eqref{defZgamma2} for $n,m\to \infty$ and the
properties $z^\gamma (\R_+)=\R_+$ and $z^\gamma (i\R_+)=\text{e}^{\gamma\pi
i/2}\R_+$
of the holomorphic mapping $z^\gamma$ motivate the following
definition of the discrete analog.

\begin{definition}\label{defZgamma}
  For $0<\gamma<2$, the discrete conformal map $Z^\gamma: \Z_+^2\to \C$
  is the solution of equations~\eqref{defZgamma1} and~\eqref{defZgamma2} with
  the initial conditions
  \begin{equation*}
    Z^\gamma(0,0)=0,\quad Z^\gamma(1,0)=1,\quad
    Z^\gamma(0,1)=\text{e}^{\gamma\pi i/2}.
  \end{equation*}
\end{definition}

From this definition, the properties $Z^\gamma(n,0)\in\R_+$ and
$Z^\gamma(0,m)\in \text{e}^{\gamma\pi i/2}\R_+$ are obvious for all
$n,m\in\N$.
Furthermore, the discrete conformal map $Z^\gamma$ from
Definition~\ref{defZgamma}
determines an $SG$-circle pattern.
Indeed, by Proposition~1 of \cite{AB00} all edges at the
vertex $f_{n,m}$ with $n+m=0\pmod{2}$ have the same length and all
angles between neighboring edges at the vertex $f_{n,m}$ with
$n+m=1\pmod{2}$ are equal to $\pi/2$. Thus, all elementary
quadrilaterals $(f_{n,m}, f_{n+1,m}, f_{n+1,m+1}, f_{n,m+1})$ build
orthogonal kites and for any $(n,m)\in \Z_+^2$ with $n+m=0\pmod{2}$
the points $f_{n+1,m}, f_{n,m+1}, f_{n-1,m}, f_{n,m-1}$ lie on a
circle with center $f_{n,m}$.
Therefore, we
consider the sublattice $\{(n,m)\in\Z_+^2: n+m=0\pmod{2}\}$ and denote by
$\V$ the quadrant
\[ \V =\{z=N+iM: N,M\in \Z, M\geq |N|\}\subset \Z^2,\]
where $ N=(n-m)/2,\quad M=(n+m)/2$. Two vertices $z_1,z_2\in\V$ are connected
by an edge if and only if $|z_1-z_2|=1$.

\begin{theorem}[{\cite{AB00,Ag03,Ag05}}]\label{theoZgamma1}
\begin{enumerate}[(i)]
 \item\label{RotDirection}
If $R(z)$ denotes the radius function corresponding to the discrete
  conformal map $Z^\gamma$ for some $0<\gamma<2$, then it holds that
  \begin{equation}
    (\gamma-1)(R(z)^2-R(z-i)R(z+1))\geq 0
  \end{equation}
  for all $z\in \V\setminus\{\pm N+iN|N\in\N\}$.
\item\label{ZgammaEinbettung}
For $0<\gamma<2$, the discrete conformal maps $Z^\gamma$ given by
Definition~\ref{defZgamma} are embedded.
Consequently, the corresponding circle patterns are also embedded.
\end{enumerate}
\end{theorem}

In the following section and in Section~\ref{SecquasiZgammaunique} we
continue to
use the notation of this section.
In particular the radius function is denoted by $R$ and
we have the normalization $R(0)=1$.

%%%%%%%%%%%%%%%%%%%%%%%%%%%%%%%%%%%%%%%%%%%%%%%%%%%%%%%%
\subsection{Uniqueness of the orthogonal $Z^\gamma$-circle
  patterns}\label{secZgammaUniq}

This section is devoted to the proof of following uniqueness result.

\begin{theorem}[Rigidity of orthogonal $Z^\gamma$-circle
patterns]\label{theoZgammauniq}
For $\gamma\in(0,2)$ the infinite orthogonal embedded circle
pattern corresponding to $Z^\gamma$ is the unique embedded orthogonal circle 
pattern (up to global scaling) with the following two properties.
\begin{enumerate}[(i)]
 \item The union of the corresponding kites of the
$Z^\gamma$-circle pattern  covers the infinite
sector $\{z=\rho\text{e}^{i\beta}\in\C : \rho\geq 0,\ \beta\in[0,\gamma\pi/2]\}$
with angle $\gamma\pi/2$.
\item The centers of the boundary circles lie on the boundary half lines $\R_+$
and $\text{e}^{i\gamma\pi/2}\R_+$.
\end{enumerate}
\end{theorem}

Our proof uses results of discrete potential theory or of the theory of
random walks which can be found in standard textbooks,
for example by Doyle and Snell~\cite{DS84} or by Woess~\cite{Woe00}.
We recall some basic terminology and notation and cite adapted versions of
a few theorems which will be useful for our argumentation.

By abuse of notation, we denote by $\Z^2$ the points $(a,b)\in\R^2\cong \C$ with
$a,b\in\Z$ as well as the graph with vertices at these points and edges
$e=[z_1,z_2]$ if $|z_1-z_2|=1$.
The meaning will be clear from the context.

Consider the network $(\Z^2,c)$ with {\em conductances} 
$c(e)>0$ and {\em resistances} $1/c(e)$ on the undirected edges $e\in E(\Z^2)$.
Then a {\em transition probability function} $p$ is given by  
\[p(z_1,z_2):=\begin{cases}c([z_1,z_2])/\bigl( \sum_{e=[z_1,z]\in
  E(\Z^2)} c(e)\bigr) &\text{if } [z_1,z_2]\in E(\Z^2) \\  0 & \text{otherwise}
\end{cases}.\] 
The stochastic process
on $\Z^2$ given by this probability function $p$ is a  reversible
random walk or a reversible
Markov chain on $\Z^2$. The {\em simple random walk} on $\Z^2$ is given by
specifying $c(e)=1$
for all edges which leads to $p(z_1,z_2)=1/4$ if $[z_1,z_2]\in E(\Z^2)$.

Denote by $p_{\text{esc}}$ the probability that a
random walk starting at any point will never return to this point. The
network ($\Z^2,c$) is called {\em recurrent} if $p_{\text{esc}}=0$ (and {\em
transient}
otherwise). Note that $p_{\text{esc}}=1/R_{\text{eff}}$, where $R_{\text{eff}}$
denotes
the effective resistance from a point to infinity.

\begin{theorem}\label{theorec}
\begin{enumerate}[(i)]
\item The simple random walk on $\Z^2$ is recurrent.
\item Let $(\Z^2,c_1)$ and $(\Z^2,c_2)$ be two networks with conductances
  $c_1(e)>0$ and $c_2(e)>0$ on the edges. If $c_2(e)\leq c_1(e)$ for
  all edges $e\in E(\Z^2)$, then the recurrence of $(\Z^2,c_1)$
  implies the recurrence of $(\Z^2,c_2)$.
\end{enumerate}
\end{theorem}
A proof can for example be found in~\cite[Chapters 5, 7, and 8]{DS84} or
\cite[Sections 1.A, 1.B, and Corollary (2.14)]{Woe00}.

A function $f:\Z^2\to \R$ is called {\em
  superharmonic} ({\em subharmonic}) with respect to the probability
function $p$ or with respect to the conductances $c$ if for every vertex
$v\in\Z^2$ we have
 $\sum_{w\in\Z^2}p(v,w)f(w)\leq f(v)$ ($\sum_{w\in\Z^2}p(v,w)f(w)\geq f(v)$).

The following proposition
shows that the
quotient of the radius functions of two orthogonal $SG$-circle patterns is
subharmonic with respect to suitably chosen conductances.
As the statement is a special case of Proposition~\ref{propSubharmGen2}
below, we omit the proof.

\begin{proposition}\label{propSubharmGen}
  Consider two orthogonal circle patterns for
  $SG(1,0)$. Denote the radii by $\rho_j$ and $r_j$ respectively,
  where $\rho_0$ and $r_0$ denote the radii of the inner circles.
Then 
  \begin{equation}
    \sum_{j=1}^4 c_j \frac{r_j}{\rho_j} \geq 
    \sum_{j=1}^4 c_j \frac{r_0}{\rho_0} \qquad \text{and}\qquad \sum_{j=1}^4 c_j
    \frac{\rho_j}{r_j} \geq  \sum_{j=1}^4 c_j \frac{\rho_0}{r_0},
  \end{equation}
where $c_j=1/((\rho_j/\rho_0)+(\rho_0/\rho_j))$.
\end{proposition}

Our proof of rigidity is based on the following property of superharmonic
functions on recurrent networks.
\begin{theorem}[{\cite[Theorem~(1.16)]{Woe00}}]\label{theosuperharm}
A network is recurrent if and only if all nonnegative superharmonic
functions are constant.
\end{theorem}

\begin{proof}[Proof of Theorem~\ref{theoZgammauniq}.]
Let $\gamma\in(0,2)$ and denote by $R:\V\to\R_+$ the radius function
of the embedded $Z^\gamma$-circle pattern
with $R(0)=1$.
Let $r:\V\to\R_+$ denote the radius function of an embedded orthogonal
circle pattern with 
the same combinatorics and the same boundary conditions (orthogonal
boundary circles to the half lines $\R_+$ and
$\text{e}^{i\gamma\pi/2}\R_+$). Without loss of generality we assume 
same normalization $r(0)=1$. This can always be achieved by a suitable scaling.
In the following, we will show that the radius functions $R$ and $r$
take the same values on all of $\V$. This implies that both circle
patterns coincide.

As both circle patterns are embedded,
Lemma~\ref{corSn} implies that for some constant $A>0$ and $n\geq 3$
\begin{equation}\label{eqrngeneration}
1-\frac{A}{n} \leq \frac{r(z_j^{(n)})}{r(z_0^{(n)})} \leq 1+\frac{A}{n}
\end{equation}
holds for all radii $r(z_0^{(n)})$ for vertices $z_0^{(n)}$ of the
$n$th generation away from the origin and 
their incident vertices $z_j^{(n)}$. Here, vertices $z\in\V$ belong to
the $n$th generation if their combinatorial distance in $\V$ to the origin is
$n$.
For estimation~\eqref{eqrngeneration} we
have also used that the reflection of the 
circle pattern in one of the boundary lines $\R_+$ or
$\text{e}^{i\gamma\pi/2}\R_+$ also leads to an embedded orthogonal $SG$-circle
pattern. The same reasoning applies to the radii of the
$Z^\gamma$-circle pattern, so
\begin{equation}\label{eqrhosn}
1-\frac{A}{n} \leq \frac{R(z_j^{(n)})}{R(z_0^{(n)})} \leq 1+\frac{A}{n}
\end{equation}
for $n\geq 3$ with the same constant $A$.
Estimations~\eqref{eqrngeneration} and~\eqref{eqrhosn}, the boundary conditions
and a suitable adaptation of the Ring Lemma~\ref{lemRingGen} for circles of
generation two and three from the origin imply that there is a constant
$K>0$ such that
\begin{equation}\label{eqrsn}
\frac{1}{K} \leq \frac{r(z_j)}{r(z_0)} \leq K\qquad \text{and} \qquad
\frac{1}{K} \leq \frac{R(z_j)}{R(z_0)}\leq K
\end{equation}
for all incident vertices $z_j$ and $z_0$.

We now consider two undirected networks $(\Z^2,C)$ and $(\Z^2,\tilde{C})$ as
follows.
On the edges of $\V$ we define two conductance functions $C$ and $\tilde{C}$ by
\begin{align*}
C(e)&=
C(R(z_j),R(z_k))
:=\left(\frac{R(z_j)}{R(z_k)} +
  \frac{R(z_k)}{R(z_j)}\right)^{-1}, \\
\tilde{C}(e) &= \tilde{C}(r(z_j),r(z_k)) :=\left(\frac{r(z_j)}{r(z_k)} +
  \frac{r(z_k)}{r(z_j)}\right)^{-1}, 
\end{align*}
where the edge $e=[z_j,z_k]$ connects the vertices $z_j,z_k\in\V$. 
Estimations~\eqref{eqrsn} imply that both positive functions
$C>0$ and
$\tilde{C}>0$ are uniformly bounded away from $0$ (and from infinity). 
These two conductance networks on $\V\subset\Z^2$ can be continued to all
of $\Z^2$ by reflection in the lines $\{|M|=|N|\}$. 
From Theorem~\ref{theorec} we deduce
that both networks $(\Z^2,C)$ and $(\Z^2,\tilde{C})$ are reccurent.

Consider the following positive functions on $\V$
\begin{equation*}
f_1(z) = r(z)/R(z)>0 \qquad \text{and}\qquad f_2(z) = R(z)/r(z)=1/f_1(z)>0.
\end{equation*}
By Proposition~\ref{propSubharmGen}
these functions are subharmonic.
Using the boundary
conditions of the circle patterns, this remains
true if $f_1$ and $f_2$ are continued to all of $\Z^2$ using
reflection. Consequently, $M-f_1$ and $M-f_2$ are superharmonic for all
constants $M\in\R$. If $f_1$ or $f_2$ is bounded from above, we
get a positive superharmonic function using the upper bound. Then
Theorem~\ref{theosuperharm} implies that
both functions are constant. Thus $r\equiv
R$ and consequently both circle patterns coincide.

Denote by $M_1(n)$ and $M_2(n)$ the maximum of $f_1$
and $f_2$, respectively, for the set of
vertices of the $n$th generation about the origin. As $f_1$ and $f_2$ are
subharmonic, they assume their maxima on the boundary. Therefore
the functions $M_1$ and $M_2$ are monotonically
increasing. The estimations~\eqref{eqrngeneration} and~\eqref{eqrhosn} imply
that the quotients of any two radii of one circle pattern in the $n$th
generation are bounded 
from above for $n\geq 3$, as two vertices in the $n$th generation can be
connected
by at most $4n$ edges using only vertices of the $n$th and $n+1$st
generation. So their quotient is bounded by $\text{e}^{4A}$ for both
radius functions $r$ and $R$. Note that with the normalization $r(0)=1=R(0)$,
the maxima $M_1$ and $M_2$ are bounded from below by $1$. Thus their product
\begin{equation*}
M_1(n)M_2(n)= f_1(z_{M_1}^{(n)})f_2(z_{M_2}^{(n)})
=\frac{R(z_{M_1}^{(n)})}{r(z_{M_1}^{(n)})}
\frac{r(z_{M_2}^{(n)})}{R(z_{M_2}^{(n)})} \leq \text{e}^{8A}
\end{equation*}
is bounded from above. Here $z_{M_1}^{(n)}$ and $z_{M_2}^{(n)}$ denote the
vertices of the $n$th generation where $f_1$ and $f_2$ assume their
maxima, respectively. Therefore $M_1$ and $M_2$ are also
bounded. This finishes the proof of uniqueness.
\end{proof}

%%%%%%%%%%%%%%%%%%%%%%%%%%%%%%%%%%%%%%%%%%%%%%%%%%%%%%%%%%%%%%%%%%%%%
\section{Uniqueness of isoradial quasicrystallic circle
patterns}\label{secQuasiCirc}

The uniqueness result of Theorem~\ref{theoSGUniq}
can be generalized for some classes of quasicrystallic circle patterns. On this
basis we will then generalize Theorem~\ref{theoZgammauniq} for some classes of
quasicrystallic $Z^\gamma$-circle pattern.

\subsection{Quasicrystallic circle patterns and connection to
 $\Z^d$}\label{secQuasiZd}

\begin{definition}
A {\em rhombic embedding in $\C$} of a b-quad-graph ${\mathscr
  D}$  is an embedding with the property that each face of
${\mathscr D}$ is mapped to a rhombus.
Given a rhombic embedding of ${\mathscr D}$, consider for each
directed edge $\vec{e}\in \vec{E}({\mathscr D})$ the vector of its embedding
as  a  complex number with length one. Half of the number of
different values of these directions is called the
{\em  dimension} $d$ of the rhombic embedding.
If $d$ is finite, the  rhombic embedding is called {\em quasicrystallic}.

Adding circles with centers in the white vertices of the rhombic embedding and
radius equal to the edge length reveals the close connection to embedded
isoradial circle patterns.

A circle pattern for a b-quad-graph
${\mathscr D}$ is called a {\em quasicrystallic circle pattern} if there
exists a quasicrystallic rhombic embedding of ${\mathscr D}$ and if the
intersection angles are taken from this rhombic embedding.
The comparison function of the isoradial circle pattern ${\mathscr 
C}_1$ for ${\mathscr D}$
and the quasicrystallic circle pattern ${\mathscr C}_2$ is also
called {\em comparison function for ${\mathscr C}_2$}.
\end{definition}

In the following we will often identify the b-quad-graph
${\mathscr D}$ with a rhombic embedding of ${\mathscr D}$.

\begin{remark}
The notion ``quasicrystallic'' is not uniquely defined in
literature. Here we adopt the definition given in~\cite{BMS05}.
Naturally, this property only makes sense for infinite graphs or
sequences of graphs with growing number of vertices and edges.
\end{remark}

\begin{floatingfigure}[l]{4.9cm}
\begin{center}
\includegraphics[height=4cm]{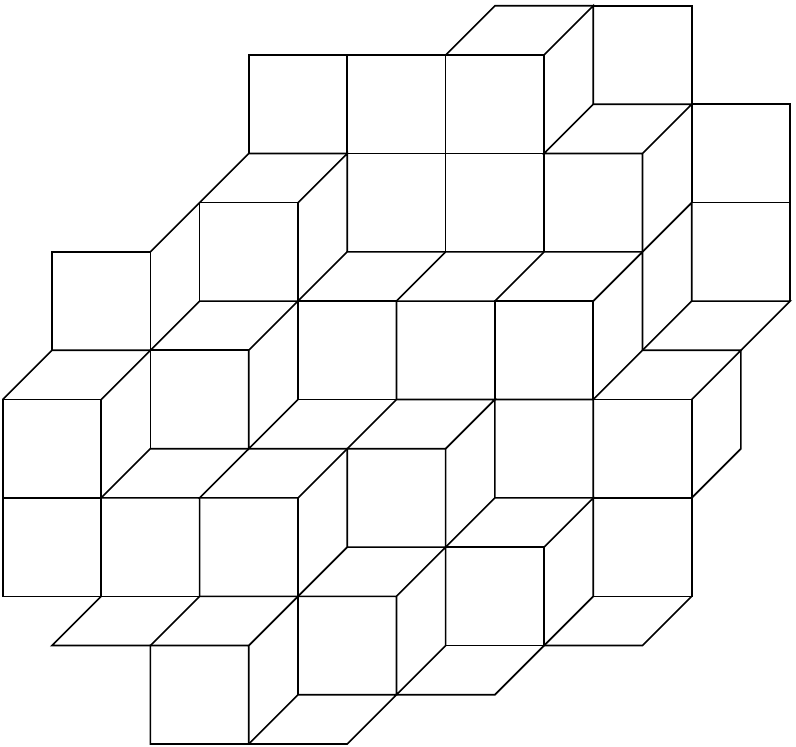}
\end{center}
\caption{An example of a quad-surface $\Omega_{\mathscr
    D}\subset \Z^3$.}\label{figExOmegaD}
 \end{floatingfigure}
Any rhombic embedding of a b-quad-graph ${\mathscr D}$ can be seen as
a sort of projection of a certain two-dimen\-sio\-nal subcomplex ({\em
quad-surface}) $\Omega_{\mathscr D}$ of the multi-dimen\-sio\-nal lattice
$\Z^d$ (or of a multi-dimen\-sio\-nal lattice $\cal L$ which is
isomorphic to $\Z^d$). An illustrating example is given in 
Figure~\ref{figExOmegaD}.

The quad-surface
$\Omega_{\mathscr D}$ in $\Z^d$ can be constructed in the following way.
Denote the set of the different edge directions of $\mathscr D$ by ${\cal
A}=\{\pm
a_1,\dots,\pm a_d\}\subset \Sp^1$. We suppose that $d>1$ and that any two
non-opposite elements of ${\cal A}$ are linearly independent over $\R$.
Let ${\mathbf e}_1,\dots,{\mathbf e}_d$ denote the standard orthonormal basis of
$\R^d$. Fix a white vertex $x_0\in V({\mathscr D})$ and
the origin of $\R^d$. Add the edges of $\{\pm {\mathbf e}_1,\dots,\pm {\mathbf
e}_d\}$ at the origin which correspond to
the edges of $\{\pm a_1,\dots,\pm a_d\}$ incident to $x_0$ in
$\mathscr D$, together with their endpoints. Successively continue the
construction at the new endpoints. Also, add two-dimensional facets
(quadrilateral faces) of
$\Z^d$ corresponding to the quadrilateral faces of $\mathscr D$, spanned by
incident edges.

A quad-surface $\Omega_{\mathscr D}$ in $\Z^d$
corresponding to a quasicrystallic rhombic embedding can be
characterized using the following monotonicity property.
For a proof see~\cite[Section~6]{BMS05}.
\begin{lemma}[Monotonicity criterium]\label{lemMonoton}
Any two points of $\Omega_{\mathscr D}$ can be connected by a path in
$\Omega_{\mathscr D}$ with all
directed edges lying in one $d$-dimensional octant, that is all
directed edges of this path are elements of one of the $2^d$ subsets of $\{\pm
{\mathbf e}_1,\dots, \pm {\mathbf e}_d\}$ containing $d$ linearly independent
vectors. 
\end{lemma}

An important class of examples of rhombic embeddings of b-quadgraphs
can be constructed using ideas of the grid projection method
for quasiperiodic tilings of the plane; see for example~\cite{DK,GR,Se}.

\begin{example}[Quasicrystallic rhombic embedding obtained from a
plane]\label{exquasi}
 Let $E$ be a two-dimensional plane in $\R^d$.
Let ${\mathbf e}_1,\dots,{\mathbf e}_d$ denote the standard orthonormal basis
of $\R^d$ and let ${\mathbf t}\in E$. We assume
that $E$ does not contain any of the segments $s_j=\{ {\mathbf t}
+\lambda{\mathbf e}_j : \lambda \in [0,1]\}$ for $j=1,\dots,d$.
If $E$ contains two different segments $s_{j_1}$ and $s_{j_2}$, the following
construction only leads to the
standard square grid pattern $\Z^2$.  If $E$ contains exactly one segment
$s_j$, the construction may be adapted for the
remaining dimensions (excluding ${\mathbf e}_j$).
We further assume that the orthogonal projections
onto $E$ of the two-dimensional facets $E_{j_1,j_2}= \{\lambda_1
{\mathbf e}_{j_1}+\lambda_2{\mathbf e}_{j_2}: \lambda_1, \lambda_2\in[0,1]\}$
for $1\leq j_1<j_2\leq d$ are non-degenerate parallelograms.
Then we can choose positive
numbers $c_1,\dots, c_d$ such that the orthogonal projections
$P_E(c_j{\mathbf e}_j)$ have length 1.

Consider around each vertex ${\mathbf p}$ of the
lattice ${\cal L}=c_1\Z\times\dots \times c_d\Z$ the hypercuboid
$V=[-c_1/2,c_1/2]\times \dots \times [-c_d/2,c_d/2]$, that is the Voronoi
cell ${\mathbf p}+V$.
These translations of $V$ cover $\R^d$.
If $E$ intersects the interior
of the Voronoi cell of a lattice point (i.e.\ $({\mathbf p}+V)^\circ\cap
E\not=
\emptyset$ for ${\mathbf p}\in {\cal L}$), then this point belongs to
$\Omega^{\cal L}(E)$. Undirected edges correspond to intersections
of $E$ with
the interior of a $(d-1)$-dimensional facet bounding two
Voronoi cells. Thus we get a connected graph
in ${\cal L}$. An intersection of $E$ with the interior
of a translated $(d-2)$-dimensional facet of $V$ corresponds to a
rectangular two-dimensional face of the lattice. By construction, the
orthogonal projection
of this graph onto $E$ results in a planar connected graph whose faces
are all of even degree ($=$ number of incident edges or of incident
vertices). A face of degree bigger than 4 corresponds to an
intersection of $E$ with the translation of a $(d-k)$-dimensional facet
of $V$ for some $k\geq 3$. Consider the vertices and edges of such a
face and the corresponding points and edges in the lattice ${\cal L}$. These
points lie on a combinatorial
$k$-dimensional hypercuboid contained in ${\cal L}$.
By construction, it is easy to see that there are two points of the
$k$-dimensional
hypercuboid which are each incident to $k$ of the given vertices. Choose
a point with least distance from $E$ and add it to the
surface. Adding edges to neighboring vertices splits
the face of degree $2k$ into $k$ faces of degree 4.

Thus we obtain an infinite monotone two-dimensional quad-surface
$\Omega^{\cal L}(E)$ which projects to an infinite rhombic embedding covering
the whole plane $E$.
Parts of such rhombic embeddings are shown in Figures~\ref{figPenrose},
\ref{figPenrose2}, and~\ref{figPenrose3}.
\end{example}

\subsection{Quasicrystallic circle patterns and 
integrability}\label{secCircInt}

Let $\mathscr D$ be a quasicrystallic rhombic embedding of a b-quad-graph.
The quad-surface $\Omega_{\mathscr D}$ in $\Z^d$ is
important by its connection with integrability. See also~\cite{BS08} for a
more detailed presentation and a deepened study of integrability and
consistency.

In particular, a
function defined on the vertices of $\Omega_{\mathscr D}$ which satisfies some
3D-consistent equation on all faces of $\Omega_{\mathscr D}$
can be uniquely extended to the {\em brick}
\begin{equation*}
\Pi(\Omega_{\mathscr D}) :=\{{\mathbf n}=(n_1,\dots, 
n_d)\in \Z^d:a_k(\Omega_{\mathscr D})\leq n_k\leq b_k(\Omega_{\mathscr D}),\
k=1,\dots, d\}.
\end{equation*}
where $a_k(\Omega_{\mathscr D})= \min_{{\mathbf n}\in V(\Omega_{\mathscr
D})}n_k$ and $b_k(\Omega_{\mathscr D})= \max_{{\mathbf n}\in V(\Omega_{\mathscr
D})}n_k$.
Note that $\Pi(\Omega_{\mathscr D})$ is the hull of
$\Omega_{\mathscr D}$. A proof may be found in~\cite[Section~6]{BMS05}.
Let $\mathscr D$ be a quasicrystallic rhombic embedding and let
${\mathscr C}_2$
be a quasicrystallic circle pattern
with the same combinatorics and the same intersection angles.
Denote the comparison function for ${\mathscr C}_2$
by $w$ as in~\eqref{eqdefw}. Since the Hirota equation~\eqref{eqw}
is 3D-consistent (see Sections~10 and~11 of~\cite{BMS05})
$w$ considered as a function on
$V(\Omega_{\mathscr D})$ can uniquely be extended to the brick
$\Pi(\Omega_{\mathscr D})$ such that equation~\eqref{eqw} 
holds on all two-dimensional facets.
Additionally, $w$ and its extension are real valued on
white points of $V(\Omega_{\mathscr D})$ and have values in $\Sp^1$ for black
points of $V(\Omega_{\mathscr D})$.
This can easily be deduced from the Hirota Equation~\eqref{eqw}.

The extension of $w$ can be used to define a radius function for any rhombic
embedding with the same boundary faces as $\mathscr D$.

\begin{lemma}\label{lem2patt}
Let ${\mathscr D}$ and ${\mathscr D}'$ be two simply connected finite
rhombic embeddings of b-quad-graphs with the same edge directions. Assume that
${\mathscr D}$ and ${\mathscr D}'$ agree on all boundary faces. Let
${\mathscr C}$ be an (embedded) planar circle pattern for ${\mathscr D}$
and the labelling given by the rhombic embedding. Then there is an
(embedded) planar circle
pattern ${\mathscr C}'$ for ${\mathscr D}'$ and the corresponding labelling
which agrees with ${\mathscr C}$ for all boundary circles.
\end{lemma}
\begin{proof}
Consider the monotone quad-surfaces $\Omega_{\mathscr D}$ and
$\Omega_{\mathscr D}'$. Without loss of generality, we can assume that
$\Omega_{\mathscr D}$ and $\Omega_{\mathscr D}'$ have the same
boundary faces in $\Z^d$. Thus both define the same brick
$\Pi(\Omega_{\mathscr D}) =\Pi(\Omega_{\mathscr D}')=:\Pi$. Given the
circle pattern ${\mathscr C}$, define the comparison function $w$ for 
${\mathscr C}$ on
$V(\Omega_{\mathscr D})$ by~\eqref{eqdefw}. Extend $w$ to the brick
$\Pi$ such that condition~\eqref{eqw} holds for all two-dimensional
facets. Consider $w$ on $\Omega_{\mathscr D}'$ and build the
corresponding pattern ${\mathscr C}'$, such that the points on the boundary
agree with
those of the given circle pattern ${\mathscr C}$. Equation~\eqref{eqw}
guarantees that all faces of $\Omega_{\mathscr D}'$
are mapped to closed kites. Due to the combinatorics, the chain of
kites is closed around each vertex. Since the boundary kites of ${\mathscr C}'$
are given by
${\mathscr C}$ which is an immersed circle pattern,
at every interior white point the angles of 
the kites sum up to $2\pi$. Thus ${\mathscr C}'$ is an immersed circle pattern.

Furthermore, ${\mathscr C}'$ is embedded if ${\mathscr C}$ is, because
${\mathscr C}'$ is an immersed circle pattern, and ${\mathscr C}'$ and
${\mathscr C}$ agree for all boundary kites.
\end{proof}

%%%%%%%%%%%%%%%%%%%%%%%%%%%%%%%%%%%%%%%%%%%%%%%%%%%%%%%%%%%%%%%%%%%%
\subsection{Local changes of rhombic embeddings}\label{secQuasiHat}

\begin{figure}[th]
 \begin{center}
 \includegraphics[height=1.5cm]{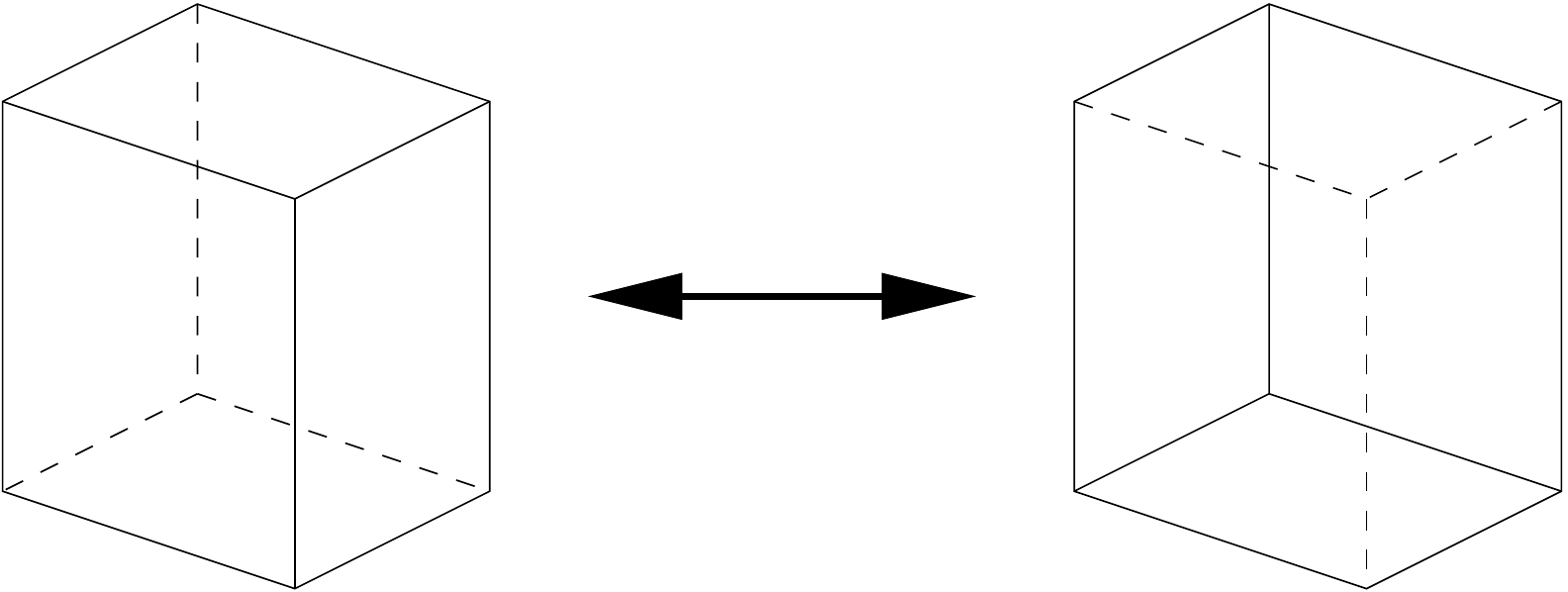}
 \end{center}
 \caption{A flip of a three-dimensional cube.
Only the faces bounded by solid edges are part of
the quad-surface in $\Z^d$.}\label{figflip}
\end{figure}

Let ${\mathscr D}$ be a rhombic embedding of a finite simply connected
b-quad-graph and let $\Omega_{\mathscr D}$ be the
corresponding quad-surface in $\Z^d$.
Let $\hat{\mathbf z}\in V_{int}(\Omega_{\mathscr D})$ be an interior vertex
with exactly three incident two-dimensional facets of $\Omega_{\mathscr D}$.
Consider the
three-dimen\-sio\-nal cube with these boundary facets. Replace 
the three given facets with the three other two-dimensional facets of
this cube. This procedure is called a {\em flip}; see Figure~\ref{figflip} for
an illustration.

A vertex ${\mathbf z}\in\Z^d$ {\em can be  reached with
flips from $\Omega_{\mathscr D}$} if ${\mathbf z}$ is contained in a
quad-surface obtained
from $\Omega_{\mathscr D}$ by a suitable sequence of flips.
The set of all
vertices which can be reached with flips, including $V(\Omega_{\mathscr D})$,
will be denoted by ${\cal F}(\Omega_{\mathscr D})$.

\begin{remark}\label{remExquasi}
The quad-surface $\Omega^{\cal L}(E)$ can also be generalized
using a finite sequence of flips.
 Such an infinite rhombic embedding will be called a
{\em plane based quasicrystallic rhombic embedding}.
\end{remark}

\begin{figure}[th]
\setlength{\unitlength}{0.75cm}
\begin{picture}(15.2,1.9)(0.6,1.2)
 \put(1,1.99){$\dots$}
\put(1,2.99){$\dots$}
\put(0.6,1.24){$\dots$}
\put(2,2){\vector(-1,0){0.5}}
\put(1.2,2.25){${\mathscr R}_+$}

\put(1.5,2){\line(1,0){5.7}}
\put(1.5,3){\line(1,0){5.7}}
\put(1.1,1.25){\line(1,0){5.7}}
\put(2,2){\line(0,1){1}}
\put(2,2){\line(-2,-3){0.49}}
\put(3,2){\line(0,1){1}}
\put(3,2){\line(-2,-3){0.49}}
\put(4,2){\line(0,1){1}}
\put(4,2){\line(-2,-3){0.49}}

\thicklines
\put(5,2){\vector(0,1){1}}
\put(5,2){\vector(-1,0){1}}
\put(5,2){\vector(-2,-3){0.49}}
\put(5.1,2.05){$\hat{\mathbf z}$}
\put(5.1,2.6){$\vec{\mathbf e}_{j_1}$}
\put(4.4,2.17){$\vec{\mathbf e}_{j_2}$}
\put(4.85,1.5){$\vec{\mathbf e}_{j_3}$}
\put(4.1,2.6){$f_1$}
\put(3.9,1.45){$f_2$}

\thinlines
\put(6,2){\line(0,1){1}}
\put(6,2){\line(-2,-3){0.49}}
\put(7,2){\line(0,1){1}}
\put(7,2){\line(-2,-3){0.49}}

 \put(7.3,1.99){$\dots$}

\thicklines
\put(7.9,2){\vector(1,0){0.7}}

\thinlines
 \put(8.9,2.24){$\dots$}
\put(9.3,2.99){$\dots$}
\put(8.9,1.24){$\dots$}

\put(9.8,3){\line(1,0){5.7}}
\put(9.4,2.25){\line(1,0){3.4}}
\put(9.4,1.25){\line(1,0){5.7}}
\put(9.8,1.25){\line(0,1){1}}
\put(10.3,3){\line(-2,-3){0.5}}
\put(10.8,1.25){\line(0,1){1}}
\put(11.3,3){\line(-2,-3){0.5}}
\put(11.8,1.25){\line(0,1){1}}
\put(12.3,3){\line(-2,-3){0.5}}
\put(12.8,1.25){\line(0,1){1}}
\put(13.3,3){\line(-2,-3){0.5}}

\put(13.4,2.05){$\hat{\mathbf z}$}

\put(13.3,2){\line(1,0){2.2}}
\put(13.3,2){\line(0,1){1}}
\put(13.3,2){\line(-2,-3){0.49}}
\put(14.3,2){\line(0,1){1}}
\put(14.3,2){\line(-2,-3){0.49}}
\put(15.3,2){\line(0,1){1}}
\put(15.3,2){\line(-2,-3){0.49}}

 \put(15.5,1.99){$\dots$}
\end{picture}
 \caption{An example of a flip for an infinite strip.}\label{figFlipinf}
\end{figure}
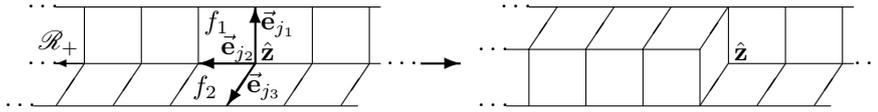

As a generalization of simple flips
we define {\em flips for simply oder doubly infinite strips} of
the
following form. See Figure~\ref{figFlipinf} for an illustration.

Let $\Omega_{\mathscr D}\subset \Z^d$ be a simply connected
monotone quad-surface. Let $\hat{\mathbf z}\in V_w(\Omega_{\mathscr
D})$ be a white vertex. Let ${\mathbf e}_{j_1}, {\mathbf e}_{j_2}, {\mathbf
e}_{j_3}$ be three different edges incident to $\hat{\mathbf z}$ such that
there are two-dimensional faces $f_1,f_2$ of $\Omega_{\mathscr D}$ incident to
${\mathbf e}_{j_1}$ and ${\mathbf e}_{j_2}$, and to ${\mathbf e}_{j_2}$ and
${\mathbf e}_{j_3}$, respectively. Let $\alpha_1=\alpha(f_1)$ and
$\alpha_2=\alpha(f_2)$ be the intersection angles associated to these faces.
Let $\alpha_3$ be the intersection angles associated to the two-dimensional
facet of $\Z^d$ incident to ${\mathbf e}_{j_1}$ and ${\mathbf e}_{j_3}$. Then
$\sum_{i=1}^3\alpha_i=2\pi$ or $\sum_{i=1}^2(\pi-\alpha_i) +\alpha_3=2\pi$. In
the first case, consider the half-axis ${\mathscr R}_+= \{\hat{\mathbf z}
+\lambda \vec{\mathbf e}_{j_2}: \lambda\geq 0\}$, where $\vec{\mathbf e}_{j_2}$
is the vector corresponding to the edge ${\mathbf e}_{j_2}$ and pointing away
from $\hat{\mathbf z}$ as in Figure~\ref{figFlipinf}. In the second case,
consider the other half-axis
${\mathscr R}_-= \{\hat{\mathbf z} +\lambda \vec{\mathbf e}_{j_2}: \lambda\leq
0\}$. In both cases we may also consider the whole axis ${\mathscr R}=
\{\hat{\mathbf z} +\lambda \vec{\mathbf e}_{j_2}: \lambda\in\R\}$. Assume that
the translations of $f_1$ and $f_2$ along these (half-)axis, that is the faces
$f_j+n\vec{\mathbf e}_{j_2} +\hat{\mathbf z}$ for $j=1,2$ and $n\in\N$,
$(-n)\in\N$, or $n\in\Z$ respectively, are contained in
$\Omega_{\mathscr D}$. We only consider the case of the positive half-axis
${\mathscr R}_+$ further. For ${\mathscr R}_-$ the argumentation is analogous
and the case of the whole axis ${\mathscr R}$ is a simple consequence. Replace
each face $f_1+n\vec{\mathbf e}_{j_2} +\hat{\mathbf z}$ by its translate
$f_1+n\vec{\mathbf e}_{j_2} +\hat{\mathbf z}+\vec{\mathbf e}_{j_3}$ for
$n\in\N_0$ and
similarly $f_2+n\vec{\mathbf e}_{j_2} +\hat{\mathbf z}$ by
$f_2+n\vec{\mathbf e}_{j_2} +\hat{\mathbf z}+\vec{\mathbf e}_{j_1}$ for
$n\in\N_0$, where $\vec{\mathbf e}_{j_1}$ and $\vec{\mathbf e}_{j_3}$
are the vectors corresponding to the edges ${\mathbf e}_{j_1}$ and ${\mathbf
e}_{j_3}$ respectively and pointing away from $\hat{\mathbf z}$. Adding the 
face incident
to $\hat{\mathbf z}$, ${\mathbf e}_{j_1}$, and ${\mathbf e}_{j_3}$, we obtain a
different, but still monotone simply connected quad-surface.

The definition of an infinite flip for a black vertex $\hat{\mathbf z}\in
V_b(\Omega_{\mathscr D})$ is very similar.

\begin{lemma}\label{lemEmbeddedInfFlip}
 Let $\Omega_{\mathscr D}\subset \Z^d$ be a simply connected
monotone quad-surface and let $\Omega_{{\mathscr D}}'$ be the
simply connected monotone quad-surface obtained from $\Omega_{\mathscr
D}$ after performing a flip for a simply infinite strip. Let $\mathscr C$ be a
circle pattern for $\mathscr D$ and the corresponding labelling and let
${\mathscr C}'$ be the corresponding circle pattern after performing the
corresponding infinite flip as for $\Omega_{\mathscr D}'$. Then the
 resulting circle pattern ${\mathscr C}'$ is embedded if the original
one $\mathscr C$ is.
\end{lemma}
The proof is based on similar arguments as the proof of Lemma~\ref{lem2patt}
and is therefore left to the reader.

%%%%%%%%%%%%%%%%%%%%%%%%%%%%%%%%%%%%%%%%%%%%%%%%%%%%%%%%%%%%%
\subsection{Uniqueness of isoradial quasicrystallic circle patterns}

Let $\mathscr E$ be the family of all infinite quasicrystallic rhombic
embeddings of connected and simply connected b-quad-graphs $\mathscr D$ which
cover the entire complex plane
and such that the brick $\Pi(\Omega_{\mathscr D})$ of the
corresponding quad-surface $\Omega_{\mathscr D}$ contains a
$\Z^2$-sublattice, that is there are at least two different indices $j_1,j_2$
such that
$\min_{\mathbf n\in V(\Omega_{\mathscr D})} n_{j_k}=-\infty$ and 
$\max_{\mathbf n\in 
V(\Omega_{\mathscr D})} n_{j_k}=\infty$ for $k=1,2$. Note that $\mathscr E$
contains in particular the plane based rhombic embeddings, like the
Penrose tilings, for which $\Pi(\Omega_{\mathscr D})=\Z^d$.
Now we use the uniqueness of $SG$-circle patterns of Theorem~\ref{theoSGUniq} 
in order to establish the uniqueness of the circle patterns of $\mathscr E$.

\begin{theorem}[Rigidity of quasicrystallic isoradial circle
patterns]\label{theoQuasiUniq}
Let $\mathscr D\in\mathscr E$ be an infinite quasicrystallic
rhombic embedding with associated graph $G$ and corresponding labelling
$\alpha$. Let $\mathscr C$ be an embedded circle pattern for $G$ and $\alpha$.
Then $\mathscr C$ is the image of the isoradial circle pattern corresponding to
$\mathscr D$ under a similarity of the complex plane.
\end{theorem}

\begin{proof}
 Let $\mathscr D\in\mathscr E$ be an infinite quasicrystallic
rhombic embeddings with associated graph $G$ and corresponding labelling
$\alpha$. Let $\mathscr C$ be an embedded circle pattern for $G$ and $\alpha$.
Consider the comparison function $w$ for $\mathscr C$ defined by~\eqref{eqdefw}
on $\Omega_{\mathscr D}$ and extend it to $\Pi(\Omega_{\mathscr D})$. Let
$\hat{\mathbf z}\in V(\Omega_{\mathscr D})$. By assumption on
$\mathscr D$, there is a $\Z^2$-sublattice $\Omega(\hat{\mathbf z})$ with
$\hat{\mathbf z}\in V(\Omega(\hat{\mathbf z}))$ which is contained in
$\Pi(\Omega_{\mathscr D})$. Furthermore, we can perform flips
for $\Omega_{\mathscr D}$ and corresponding flips for the circle pattern
$\mathscr C$ such that the
resulting quad-surface $\Omega'$ contains an arbitrary number of
generations of the lattice $\Omega(\hat{\mathbf z})$ about $\hat{\mathbf z}$.
As the corresponding circle pattern $\mathscr C'$ is embedded by
Lemma~\ref{lemEmbeddedInfFlip} and as the number of
generations about $z$ can be chosen arbitrarily large, we deduce from the
Rigidity Theorem~\ref{theoSGUniq} that the radius function is constant on
$\Omega(\hat{\mathbf z})$. More precisely, the extension of $w$ is constant on
white and black vertices of $\Omega(\hat{\mathbf z})$ respectively.

This argumentation is valid for all vertices $\hat{\mathbf z}\in
V(\Omega_{\mathscr D})$ and all $\Z^2$-sublattice $\Omega(\hat{\mathbf z})$ with
$\hat{\mathbf z}\in V(\Omega(\hat{\mathbf z}))$ which are contained in
$\Pi(\Omega_{\mathscr D})$. Therefore the extension of $w$ is constant, on
white and black vertices respectively, on all $\Z^2$-sublattices which are
contained in $\Pi(\Omega_{\mathscr D})$. Due to our assumptions on the
combinatorics of $\mathscr
D$, the radius function which is $w$ restricted to white vertices has to be
constant on the whole brick $\Pi(\Omega_{\mathscr D})$. This implies in
particular that
$\mathscr C$ is an isoradial circle pattern and thus is the image of the
isoradial circle pattern corresponding to
$\mathscr D$ under a similarity of the complex plane.
\end{proof}

\begin{definition}
Let ${\mathscr D}$ be a rhombic embedding of a simply connected
b-quad-graph.
The {\em combinatorial $K$-environment} of a point $z\in V({\mathscr D})$ is
the subgraph corresponding to all vertices which have combinatorial distance at
most $K$ from $z$ in ${\mathscr D}$.
\end{definition}

\begin{corollary}
 Let $\mathscr D\in\mathscr E$ be an infinite quasicrystallic
rhombic embeddings with corresponding labelling $\alpha$. Let $v_0\in
V_w({\mathscr D})$ be a white vertex. Then there are a constant
$n_0=n_0({\mathscr D})\in\N$ and a sequence
$s_n(v_0,{\mathscr D})$ decreasing to $0$ for $n\to\infty$ such that the
following holds.

For $n\in\N$, $n\geq n_0$, let ${\mathscr D}_{2n}(v_0)$ be the rhombic embedding
corresponding to the $2n$-environment of $v_0$. Let $G_{n}(v_0)$ be the
associated graph. Let ${\mathscr C}_{n}$ be an
embedded circle pattern for $G_{n}(v_0)$ and the labelling $\alpha$ taken from
$\mathscr D$ with radius function $r_n$.
Then there holds
\begin{equation}
 \left|\frac{r_n(v_0)}{r_n(v_1)} -1\right| \leq s_n(v_0,{\mathscr D})
\end{equation}
for all vertices $v_1\in V(G_{n}(v_0))$ incident to $v_0$.
\end{corollary}
We omit the proof which is very similar to the proof of the Hexagonal Packing
Lemma of~\cite{RS87}.
The following lemma is similar to
Lemma~\ref{lemRinggen} and corresponds to the Ring Lemma of~\cite{RS87}.

\begin{lemma}\label{lemRingGen}
 Let $\mathscr D$ be a quasicrystallic
rhombic embeddings with associated graph $G$ and
labelling $\alpha$. Denote by $\alpha_{min}=\min\{\alpha(e) : e\in E(G)\}$ the
smallest intersection angle. Let $n_0\in\N$ be such that
$(n_0-3)\alpha_{min}>\pi$. Let $v_0\in
V_w({\mathscr D})$ be a white vertex. Assume that $\mathscr D$ contains a
$(2n_0)$-environment about $v_0$. Then there is a constant $C=C(\mathscr D)>0$
such that the following holds.

Let $r$  be the radius function of an embedded circle pattern for $\mathscr D$ 
and $\alpha$ and let $v_1$ be a vertex incident to $v_0$ in $G$. Then
\[\frac{r(v_1)}{r(v_0)}>C.\] 
\end{lemma}
 \begin{proof}
Suppose that there is a vertex $v_1$ incident to $v_0$ and a sequence of
embedded circle patterns for
$\mathscr D$ and $\alpha$ with radius functions $r_n$ such that $r_n(v_0)=1$
and $r_n(v_1)\to 0$ as $n\to\infty$. Without loss of generality
we may assume that the circle $C_0$ corresponding to the vertex $v_0$ and
the intersection point corresponding to one fixed black vertex $w_0$ incident
to $v_0$ and $v_1$ in $\mathscr D$ are fixed for the whole
sequence. Then there is a subsequence such that all the circles converge to
circles or lines, that is converge in the Riemann sphere $\hat{\C}\cong
\Sp^2$.

Equation~\eqref{eqFgen} implies that there are at least two circles
corresponding to vertices incident to $v_0$ whose radii do not converge to
$0$.
If the limit is finite, there are at least two circles whose radii do not
converge to $0$ corresponding to vertices incident to this
vertex of the first generation and so on.

Consider the kites which contain in the limit the
intersection point corresponding to $w_0$ and apply the above
argument at most $n_0$ times. Then by our assumption on $n_0$ we obtain two
kites whose interiors intersect in the limit. This is a contradition to the
embeddedness of the sequence.
\end{proof}

If $\mathscr D$ is a plane based quasicrystallic rhombic embedding
we also obtain an analog of the Rodin-Sullivan Conjecture,
see~\cite{RS87,He91,Ah94}.

\begin{corollary}\label{corSnquasi}
 Let $\mathscr D$ be a plane based quasicrystallic rhombic embedding.
There are absolute constants $C=C({\mathscr D})>0$ and $n_0=n_0({\mathscr
D})\in\N$, depending only on $\mathscr D$, such that for all white vertices
$v_0\in V_w({\mathscr D})$ and all $n\geq n_0$ there holds
  \begin{equation}
    s_n(v_0,{\mathscr D}) \leq s_n({\mathscr D})\leq C/n.
  \end{equation}
\end{corollary}
\begin{proof}
 Let  $v_0\in V_w({\mathscr D})$ be any white vertex. If $n\geq
n_0({\mathscr D})$ is big enough, for each $\Z^2$-sublattice
$\Omega(\hat{\mathbf v}_0)$ the set ${\cal F}(\Omega_{{\mathscr D}})$
contains a $\lfloor B({\mathscr D})n\rfloor$-environment of $\hat{\mathbf v}_0$
in $\Omega(\hat{\mathbf v}_0)$, where the constant $B({\mathscr D})$ depends
only on the construction parameters of $\mathscr D$. 
Here $\lfloor p \rfloor$ denotes the largest integer smaller than $p\in\R$.
Therefore we can choose $s_n({\mathscr D})$ to be the maximum of $s_{\lfloor
B({\mathscr D})n\rfloor}({\mathscr D}(\hat{\mathbf v}_0))$ for all
possible regular rhombic embeddings ${\mathscr D}(\hat{\mathbf v}_0)$
corresponding to $\Z^2$-sublattice $\Omega(\hat{\mathbf v}_0)$. Now the
claim follows from Corollary~\ref{corSn} for $SG$-circle patterns.
\end{proof}

\section{Uniqueness of quasicrystallic $Z^\gamma$-circle
patterns}\label{SecquasiZgammaunique}

In this section we consider quasicrystallic $Z^\gamma$-circle patterns as
defined in~\cite{BMS05} and then prove their rigidity. The proofs are based on
the results of the previous sections and on similar arguments as for orthogonal
$Z^\gamma$-circle patterns.

\subsection{Definition and useful properties}

Let $\psi\in (0,\pi)$ be a fixed angle.
Recall the definition of the labelling $\alpha_\psi$ in~\eqref{defalphapsi}.
We consider the following generalization of
Definition~\ref{defZgamma}.

\begin{definition}[{\cite{Agpsi}}]\label{defZgammapsi}
 For $0<\gamma<2$, the discrete map $Z^\gamma: \Z_+^2\to \C$
  is the solution of
\begin{multline}
q(f_{n,m}, f_{n+1,m}, f_{n+1,m+1}, f_{n,m+1}):=  \\
\frac{(f_{n,m}- f_{n+1,m})(f_{n+1,m+1}- f_{n,m+1})}{(f_{n+1,m}-
  f_{n+1,m+1})(f_{n,m+1}- f_{n,m})} = \text{e}^{2i(\psi-\pi)}
\end{multline}
and~\eqref{defZgamma2} with  the initial conditions
  \begin{equation*}
    Z^\gamma(0,0)=0,\quad Z^\gamma(1,0)=1,\quad
    Z^\gamma(0,1)=\text{e}^{\gamma(\pi-\psi) i}.
  \end{equation*}
\end{definition}

As in the orthogonal case one can again associate a
circle pattern (for a quadrant of $SG$ corresponding to $\Z_+^2$ and related 
to $\V$ and $\alpha_\psi$) to the map $Z^\gamma$, see~\cite{Agpsi} for more
details.
Furthermore, the following results generalize the orthogonal case.

\begin{theorem}[{\cite{Agpsi}}]\label{theoZgammapsi}
\begin{enumerate}[(i)]
\item\label{RotDirectionpsi}
  If $R(z)$ denotes the radius function corresponding to the discrete
  conformal map $Z^\gamma$ for some $0<\gamma<2$, then
  \begin{equation}\label{eqRotDirectionpsi}
    (\gamma-1)(R(z)^2-R(z-i)R(z+1) -\cos\psi R(z)(R(z-i)-R(z+1))\geq 0
  \end{equation}
  for all $z\in \V\setminus\{\pm N+iN|N\in\N\}$.
\item\label{ZgammaEinbettungpsi}
For $0<\gamma<2$, the discrete conformal maps $Z^\gamma$ given by
Definition~\ref{defZgamma} are embedded.
Consequently, the corresponding circle patterns
are also embedded.
\end{enumerate}
\end{theorem}

Definition~\ref{defZgammapsi} can be generalized further.
As explained in Section~13 of~\cite{BMS05}, discrete analogs of the power
function $z^\gamma$ can also be defined for quasicrystallic rhombic embeddings
$\mathscr D$
instead of $\Z^2_+$ (or $\Z^2$). In particular, let ${\cal A}=\{\pm
a_1,\dots,\pm a_d\}\subset \Sp^1$ be the set of edge directions. Suppose that
$d>1$ and that any two
non-opposite elements of ${\cal A}$ are linearly independent over $\R$. For
$0<\gamma<2$ define the following values of the comparison function $w$ on the
coordinate semi-axis of $\Z^d_+$:
\begin{equation}\label{eqdefsemi}
 w(ne_k)=\begin{cases}
        1 &\text{if } n=0,\\
	a_k^{\gamma-1} =\text{e}^{(\gamma-1)\log a_k} &\text{if } n
\text{ is odd},\\
	\prod_{m=1}^{n/2} \frac{m-1+\frac{\gamma}{2}}{m- \frac{\gamma}{2}}
&\text{if } n\geq 2 \text{ and } n \text{ is even}.
         \end{cases}
\end{equation}
The value of the logarithm $\log a_k$ is chosen as follows.
Without loss of generality, we assume a circular order of the points of
$\cal A$ on the positively oriented unit circle $\Sp^1$ is $a_1,\dots,
a_d,-a_1,\dots, -a_d$. Set $a_{k+d}=-a_k$ for $k=1,\dots,d$ and define
$a_m$ for all $m\in\Z$ by $2d$-periodicity. To each
$a_m=\text{e}^{i\theta_m}\in\Sp^1$ assign a certain value of the argument
$\theta_m\in\R$: choose $\theta_1$ arbitrarily and then use the rule
\[\theta_{m+1}-\theta_m\in(0,\pi)\qquad \text{for all }m\in\Z. \]
Clearly we then have $\theta_{m+d}=\theta_m+\pi$.
The points $a_m$ supplied with the arguments $\theta_m$ can be considered as
belonging to the Riemann surface of the logarithmic function (i.e.\ a branched
covering of the complex plane). Now, the branch of the logarithm
is chosen such that
\[ \log(a_l)\in[i\theta_m,i\theta_{m+d-1}],\qquad l=m,\dots,m+d-1.\]

Using the Hirota Equation~\eqref{eqw}, this function $w$ can be extended to the
whole sector $\Z^d_+$. Using suitable branches of the logarithm, $w$ may also
be extended to other sectors or to a branched covering of $\Z^d$.

Figure~\ref{figVeronika} shows an example of such a quasicrystallic
$Z^\gamma$-circle pattern.

\begin{figure}[t]
\begin{center}
\includegraphics[height=
0.4\textwidth]{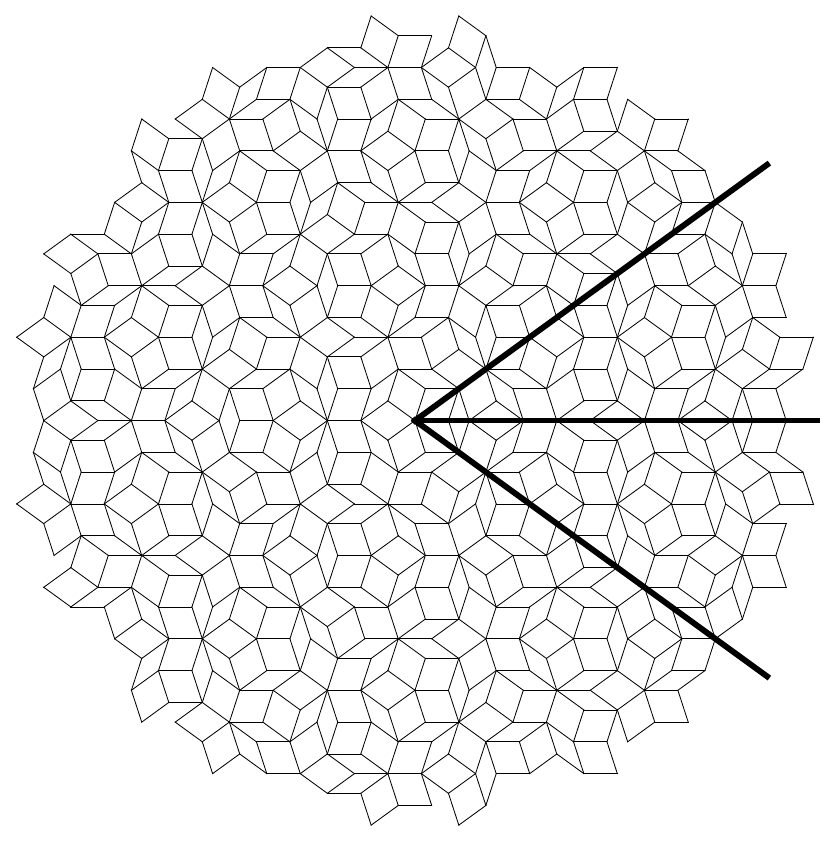}
\end{center}
\caption{A quasicrystallic rhombic embedding with
five-fold rotation symmetry.}\label{figPenrose}
\end{figure}

\begin{figure}[t]
\begin{center}
\includegraphics[height=
0.45\textwidth]{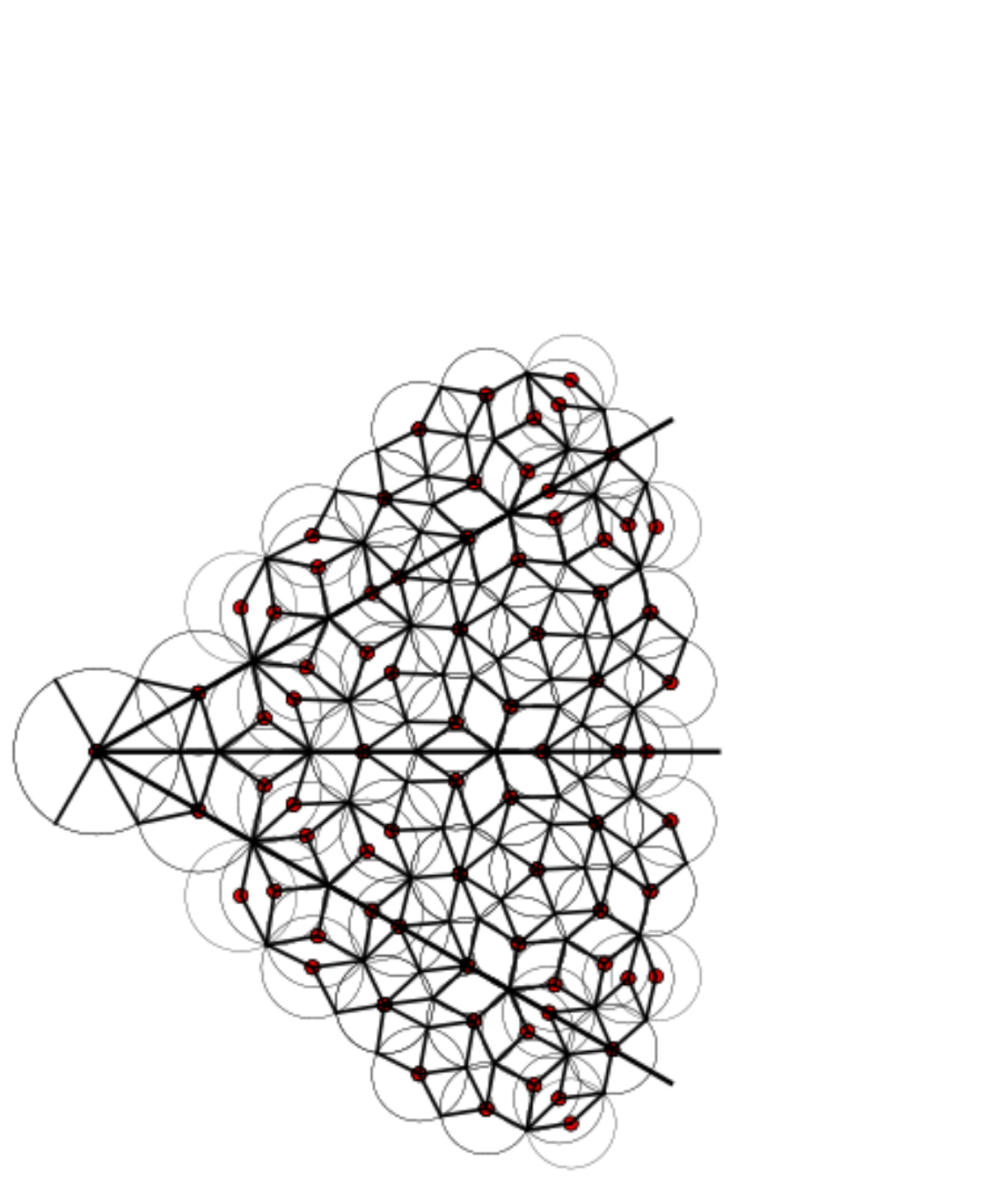}
\hspace{2em}
\includegraphics[height=
0.45\textwidth]{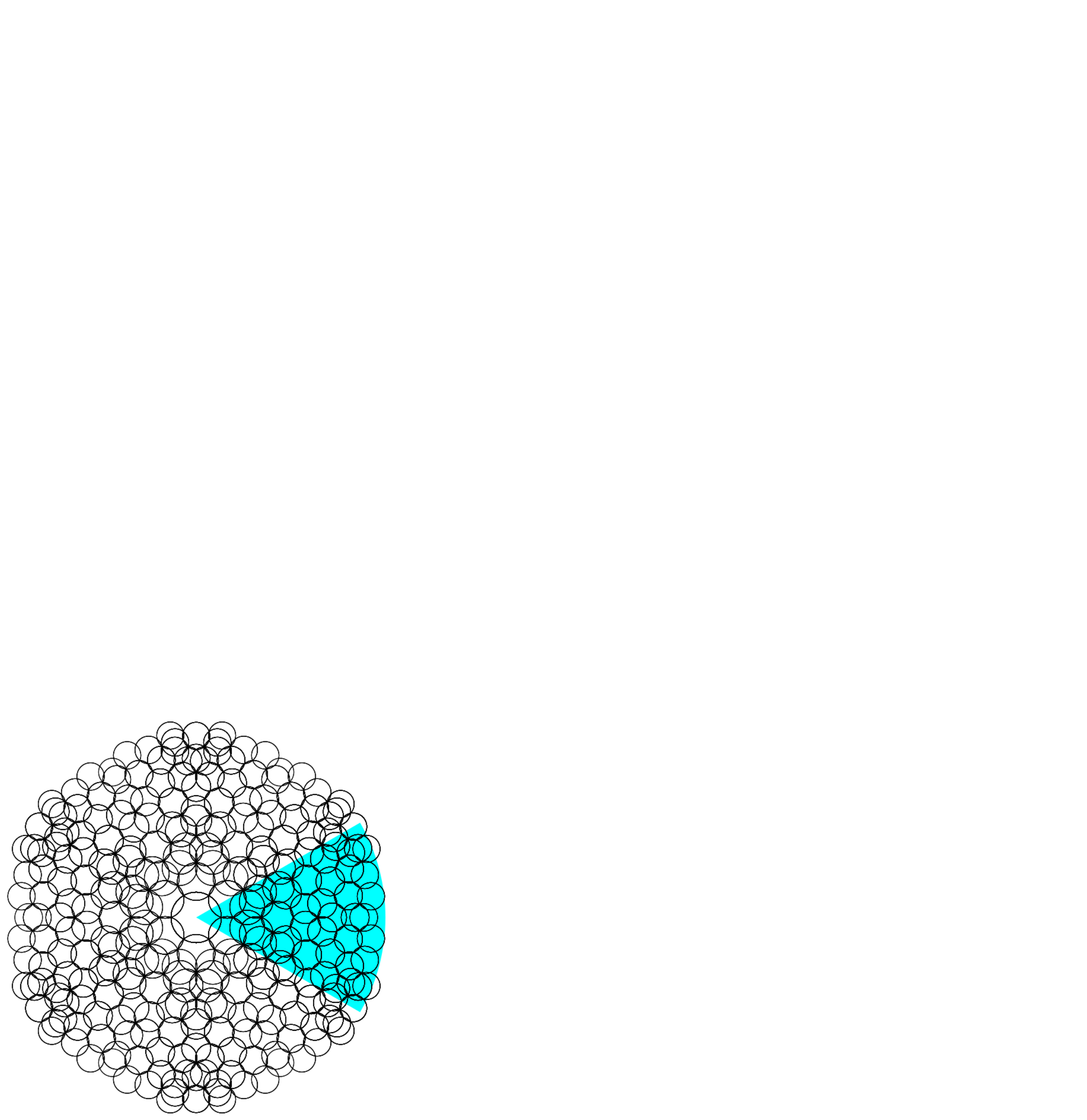}
\end{center}
\caption{Example of a quasicrystallic $Z^{5/6}$-circle pattern from the rhombic
embedding of Figure~\protect\ref{figPenrose}: Construction ({\it left}) and
corresponding circle pattern ({\it right}).}\label{figVeronika}
\end{figure}

Note that for $d=2$ the boundary conditions given in~\eqref{eqdefsemi} lead to
the circle patterns specified in Definition~\ref{defZgammapsi}.
Therefore, by Theorem~\ref{theoZgammapsi}~\eqref{ZgammaEinbettungpsi}, the
circle
pattens corresponding to the restriction of $w$ to quad-surfaces
$\Z^2_+\subset\Z^d_+$ which are spanned by two coordinate semi-axis are
embedded. We can apply finite and infinite flips to obtain other monotone
quad-surfaces corresponding to rhombic embeddings. In particular, we
obtain restrictions to $\Z^d_+$ of the plane based quad-surfaces
constructed in
Example~\ref{exquasi} and Remark~\ref{remExquasi}.
Lemmas~\ref{lem2patt} and~\ref{lemEmbeddedInfFlip} imply that these
lead again to embedded circle patterns. Thus we have proven

\begin{theorem}[Embeddedness of quasicrystallic $Z^\gamma$-circle
patterns]
Let $\Omega\subset\Z_+^d$ be a simply connected monotone quad-surface.
Then the circle pattern given by the function $w$ with initial
values~\eqref{eqdefsemi} is embedded.
\end{theorem}

\begin{figure}[tb]
\begin{center}
\includegraphics[height=
0.4\textwidth]{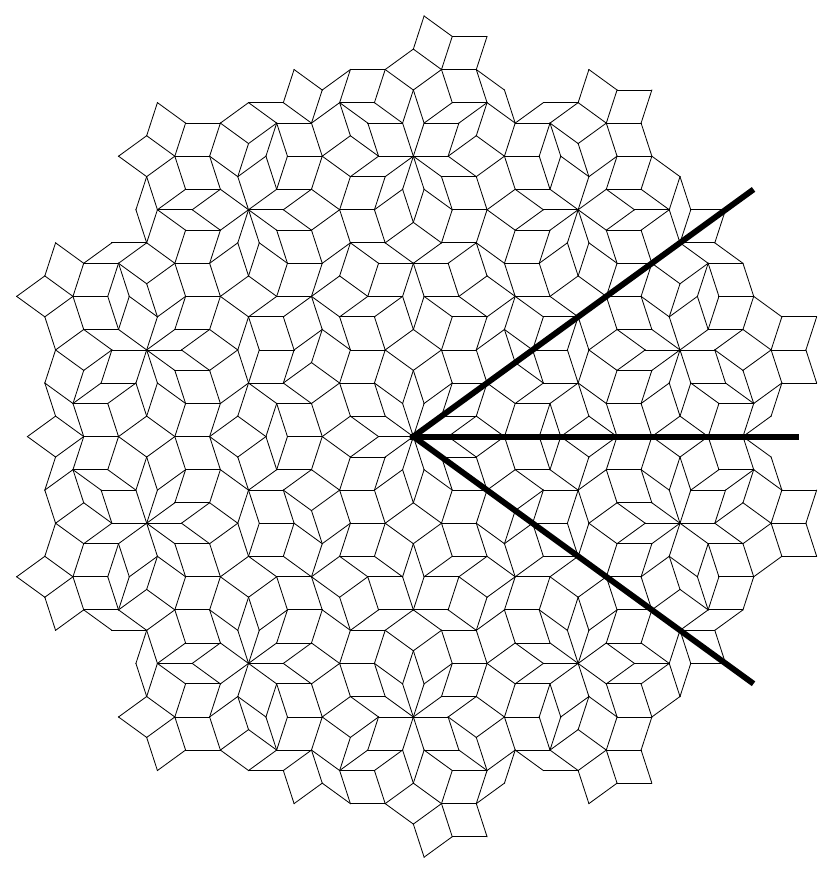}
\end{center}
\caption{A quasicrystallic rhombic embedding with
ten-fold rotation symmetry.}\label{figPenrose2}
\end{figure}

\begin{figure}[tb]
\begin{center}
\includegraphics[height=
0.5\textwidth]{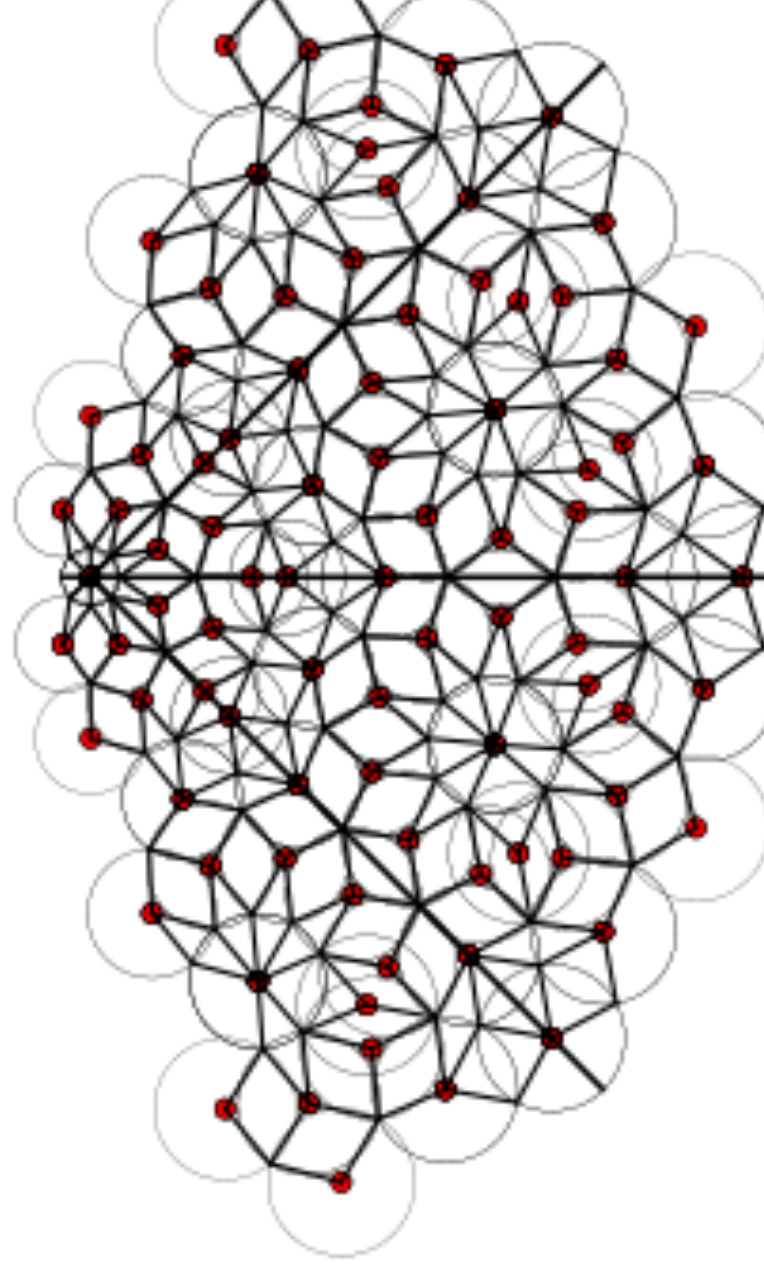}
\hspace{2em}
\includegraphics[height=
0.5\textwidth]{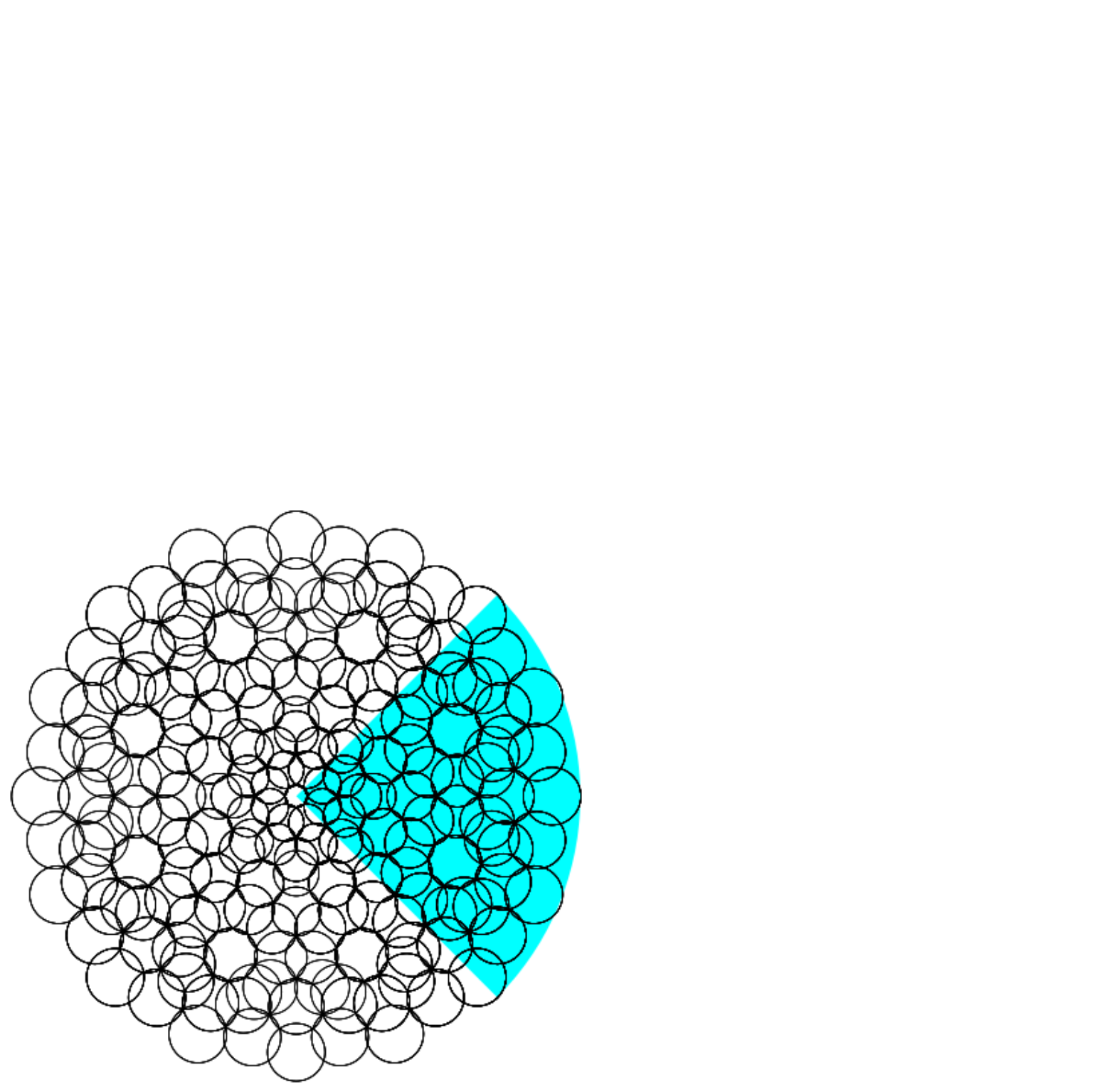}
\end{center}
\caption{Example of a quasicrystallic $Z^{5/4}$-circle pattern  from the rhombic
embedding of Figure~\protect\ref{figPenrose}: Construction ({\it left}) and 
corresponding circle pattern ({\it right}).}\label{figVeronika2}
\end{figure}

\begin{example}[Construction of the examples in Figures~\ref{figVeronika},
\ref{figVeronika2}, and~\ref{figVeronika3}]
The pictures in Figures~\ref{figVeronika}, \ref{figVeronika2}
and~\ref{figVeronika3} have been obtained using a computer program
implemented by Veronika Schreiber for her master thesis~\cite{Veronika}.

The rhombic embeddings in Figures~\ref{figPenrose}
and~\ref{figPenrose2} can be obtained by the construction method explained in
Example~\ref{exquasi} using the affine planes
$E_1=\{x=t_1+\lambda_1 u_1 +\lambda_2 u_2, \lambda_1,\lambda_2\in\R\}$ and
$E_2=\{x=t_2+\lambda_1 u_1 +\lambda_2 u_2, \lambda_1,\lambda_2\in\R\}$
respectively, where $u_1,u_2,t_1,t_2$ are defined as follows.
$\{u_1,u_2,u_3,u_4,u_5\}$ is an orthonormal basis such that the matrix
$\left(\begin{smallmatrix} 0 &0 &0 &0 &1 \\ 1&0 &0 &0 &0 \\ 0& 1&0 &0 &0 \\0 &
0& 1&0 &0 \\0&0 &0 &1& 0 \end{smallmatrix}\right)$ takes the form
\[\begin{pmatrix} \cos(2\pi/5) & -\sin(2\pi/5) & 0 &0 & 0 \\
   \sin(2\pi/5)  & \cos(2\pi/5)& 0 &0 & 0 \\
0 & 0 & \cos(4\pi/5) & -\sin(4\pi/5) & 0 \\
0 & 0 &\sin(4\pi/5)  & \cos(4\pi/5)& 0 \\
0 &0 &0 &0 &1 \end{pmatrix} \]
with respect to this basis. The translation vectors are
$t_1=\left(\begin{smallmatrix} -0.2 &-0.2 &-0.2& -0.2 &-0.2
\end{smallmatrix}\right)^T$ and $t_2=\left(\begin{smallmatrix} -0.5 &-0.5 &-0.5&
-0.5 &-0.5 \end{smallmatrix}\right)^T$ respectively. Here $v^T$ is the
transpose of the vector $v$. Translations in direction of the vector
$\left(\begin{smallmatrix} 1 & 1 & 1 & 1 & 1 \end{smallmatrix}\right)^T$
generally lead to rotationally symmetric rhombic embeddings.

For the remaining
construction it is important that the rhombic embeddings have rotational
symmetry. Now consider the two sectors indicated by lines in
Figures~\ref{figPenrose}
and~\ref{figPenrose2} or more precisely the corresponding 5-dimensional octant
in $\Z^5$. Choose an exponent $\gamma$ of the power function $z^\gamma$. In
order to construct an image circle pattern which closes up as in
% Figures~\ref{figZ5_6} and \ref{figZ5_4} %and~\ref{figZ5_7}
Figures~\ref{figVeronika}~({right}) and~\ref{figVeronika2}~({right})
we need to choose an integer $p\geq 3$ and take $\gamma=5/p$.
Given these ingredients, define the values of the comparison function $w$ on
the coordinate axes of the octant in $\Z^5$ according to~\eqref{eqdefsemi} and
calculate the missing values for the octant using the Hirota
equation~\eqref{eqw}. Taking the values on the quad-surface
corresponding to the original rhombic embedding, we can construct a sector of
the desired quasicrystallic circle pattern, see
Figures~\ref{figVeronika}~({left}) and~\ref{figVeronika2}~({left}).
The closed circle patterns in Figures~\ref{figVeronika}~({right})
and~\ref{figVeronika2}~({right}) are obtained using the rotational symmetry.

\begin{figure}[tb]
\begin{center}
\includegraphics[height=
0.55\textwidth]{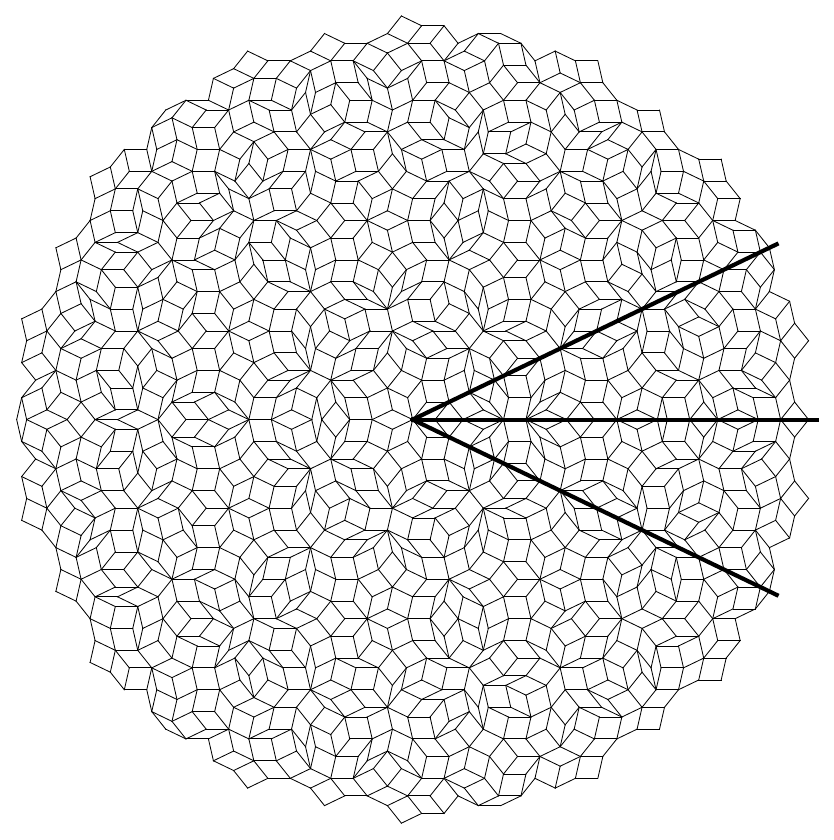}
\end{center}
\caption{A quasicrystallic rhombic embedding with
seven-fold rotation symmetry.}\label{figPenrose3}
\end{figure}

\begin{figure}[tb]
\begin{center}
\includegraphics[height=
0.55\textwidth]{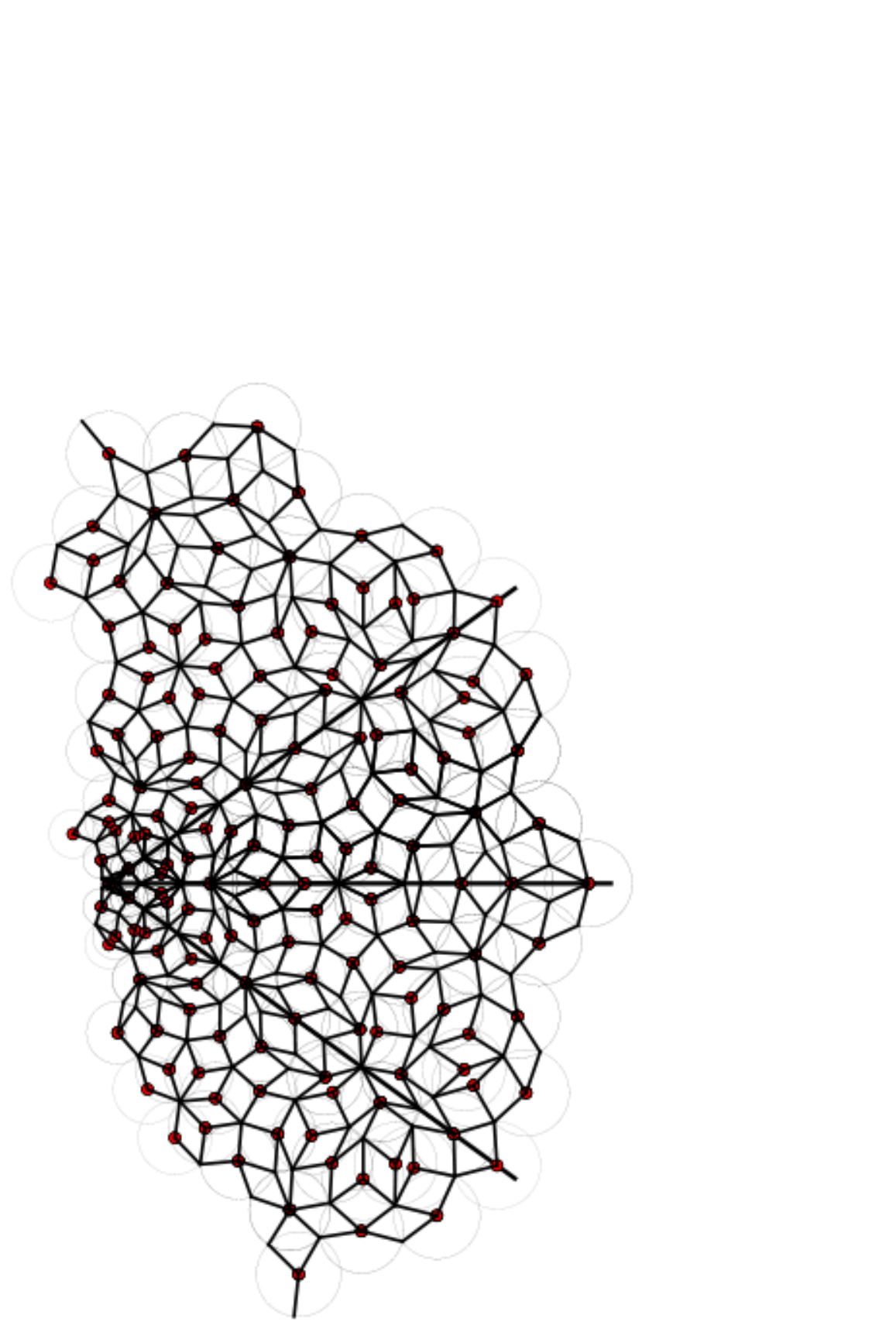}
\includegraphics[height=
0.55\textwidth]{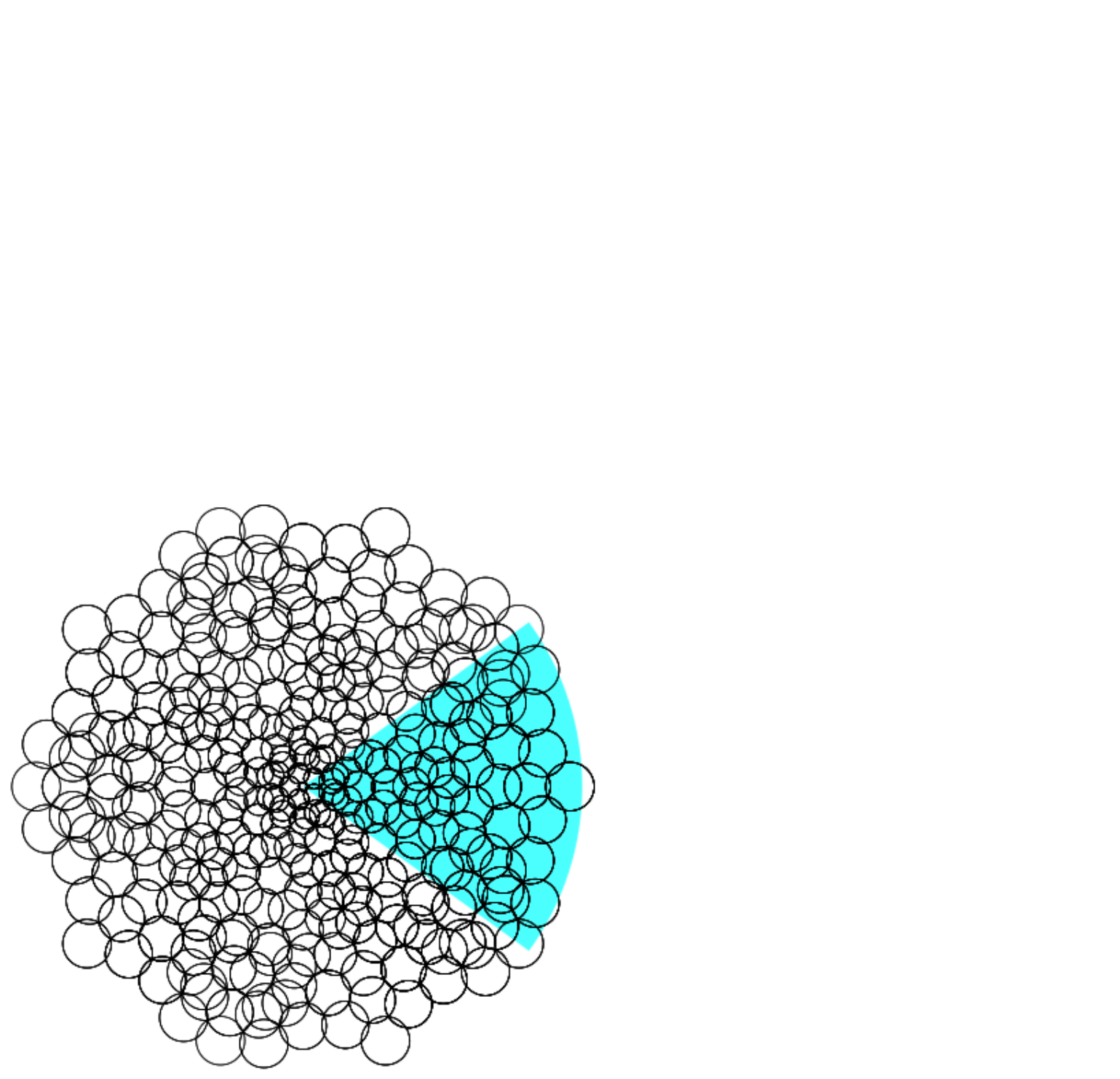}
\end{center}
\caption{Example of a quasicrystallic $Z^{7/5}$-circle pattern from the rhombic
embedding of Figure~\protect\ref{figPenrose}: Construction ({\it left}) and 
corresponding circle pattern ({\it right}).}\label{figVeronika3}
\end{figure}

The construction of the examples in Figures~\ref{figPenrose3}
and~\ref{figVeronika3} is similar using a suitable affine plane in $\Z^7$.
\end{example}

%%%%%%%%%%%%%%%%%%%%%%%%%%%%%%%%%%%%%%%%%%%%%%%%%%%%%%%%%%%%%%%%%%%
\subsection{Uniqueness of quasicrystallic $Z^\gamma$-circle
patterns}\label{secZgammaquasiUniq}

In this section we prove uniqueness of quasicrystallic
$Z^\gamma$-circle patterns using the same arguments as for the orthogonal case.

We begin with a generalization of Proposition~\ref{propSubharmGen}.
Note that we need a (geometric) restriction which
ensures that the kites corresponding to intersecting circles are convex. 
Unfortunately, flips (finite or infinite) may destroy convexity.

\begin{proposition}\label{propSubharmGen2}
Let $\mathscr D$ be a quasicrystallic rhombic embedding with associated graph
$G$. Let $\alpha$ be the labelling corresponding to $\mathscr D$. Let $v_0$ be
an interior
vertex of $G$ with incident vertices $v_1,\dots,v_m$.
  Consider two circle patterns for $G$ and $\alpha$ with radius functions $r$
and $\rho$. Denote $r_j=r(v_j)$ and $\rho_j=\rho(v_j)$ for $j=0,1,\dots,m$
and suppose that $r_j\geq r_0\cos\alpha_j$ and $\rho_j\geq \rho_0\cos\alpha_j$
for $j=1,\dots,m$, where $\alpha_j=\alpha([v_0,v_j])$. Then
  \begin{equation}
    \sum_{j=1}^4 c_j \frac{r_j}{\rho_j} \geq 
    \sum_{j=1}^4 c_j \frac{r_0}{\rho_0} \qquad \text{and}\qquad \sum_{j=1}^4 c_j
    \frac{\rho_j}{r_j} \geq \sum_{j=1}^4 c_j \frac{\rho_0}{r_0},
  \end{equation}
where $c_j=\sin\alpha_j/((\rho_j/\rho_0)+(\rho_0/\rho_j) -2\cos\alpha_j)$ for
$j=1,\dots,m$.
\end{proposition}
\begin{proof}
The proof is
based on a Taylor expansion of $f_\alpha (x+\log y)$ about $y=1$.
\[f_\alpha(x+\log y)=  f_\alpha(x) + f_\alpha'(x)(y-1)
-\frac{\sin\alpha (\text{e}^{x+\log\xi}
-\cos\alpha)}{2\xi^2(\cosh(x+\log\xi)- \cos\alpha)^2}(y-1)^2\]
by Lemma~\ref{lemPropf}~(i) with $\xi=t +(1-t)y$ for some suitable $t\in(0,1)$
depending on $x$ and $y$.
Equation~\eqref{eqFgen} for the two circle patterns implies
\begin{align*}
 \pi & = \sum_{j=1}^m f_{\alpha_j}(\log\textstyle \frac{r_j}{r_0})
= \sum_{j=1}^m f_{\alpha_j}(\log\textstyle \frac{\rho_j}{\rho_0} +\log
\textstyle \frac{r_j\rho_0}{r_0\rho_j}) \\
&= \underbrace{\sum_{j=1}^m f_{\alpha_j}(\log{\textstyle
\frac{\rho_j}{\rho_0}})}_{=\pi}
+\sum_{j=1}^m f_{\alpha_j}'(\log{\textstyle \frac{\rho_j}{\rho_0}})
\left( \frac{r_j\rho_0}{r_0\rho_j} -1\right) \\
&\hspace{8em}
-\sum_{j=1}^m\frac{(\frac{\rho_j}{\rho_0}\xi_j
-\cos\alpha_j)\sin\alpha_j}{2\xi_j^2(\cosh(\rho_j/\rho_0+\log\xi_j)-
\cos\alpha_j)^2}\left( \frac{r_j\rho_0}{r_0\rho_j} -1\right)^2,
\end{align*}
where $\xi_j=t_j +(1-t_j)\frac{r_j\rho_0}{r_0\rho_j}>0$ with
suitable $t_j\in(0,1)$ for $j=1,\dots,m$. Furthermore
\[\frac{\rho_j}{\rho_0}\xi_j -\cos\alpha_j =t_j \frac{\rho_j}{\rho_0}
+(1-t_j)\frac{r_j}{r_0}-\cos\alpha_j\geq 0.\]
by our assumption. Thus
\[\sum_{j=1}^m f_{\alpha_j}'(\log{\textstyle \frac{\rho_j}{\rho_0}})
\left( \frac{r_j\rho_0}{r_0\rho_j} -1\right)\geq 0.\]
This implies the first claim since $f_{\alpha_j}'(\log \frac{\rho_j}{\rho_0})
=c_j$ by Lemma~\ref{lemPropf}~(i).

The second claim follows from the fact,
that $1/\rho$ and $1/r$ are also radius functions of circle patterns for $G$
and $\alpha$ by Lemma~\ref{lemInvRad}.
Also, the coefficients $c_j$ are invariant
under the transformation $\rho\mapsto 1/\rho$.
\end{proof}

The following lemma specifies for which parameters
the kites of the $Z^\gamma$-circle patterns are convex. For the cases
excluded below, there exist non-convex kites, because already the kite
built from the circle centered at the origin is non-convex.

\begin{lemma}\label{lemZgammaconvex}
 If $\psi\geq \pi/2$ and $0<\gamma<2$ or if $\psi< \pi/2$ and
$(\pi-2\psi)/(\pi-\psi)\leq \gamma
\leq \pi/(\pi-\psi)$, all kites in the $Z^\gamma$-circle
pattern given by Definition~\ref{defZgammapsi} are convex.
\end{lemma}
The proof is technical. A brief version is presented in the appendix~\ref{app1}.
More details can be found in~\cite{diss}.

\begin{theorem}[Rigidity of $Z^\gamma$-circle
patterns from Definition~\ref{defZgamma1}]\label{theoZgammauniqpsi}
 If $\psi\geq \pi/2$ and $0<\gamma<2$ or if $\psi< \pi/2$ and
$(\pi-2\psi)/(\pi-\psi)\leq \gamma \leq \pi/(\pi-\psi)$ then the
$Z^\gamma$-circle pattern given by Definition~\ref{defZgammapsi} is
the unique embedded $SG$-circle
pattern for $\Z^2_+$ and $\alpha_\psi$ (up to global scaling) with the following
properties.
\begin{enumerate}[(i)]
\item The infinite sector $\{z=\rho\text{e}^{i\beta}\in\C : \rho\geq 0,\
\beta\in[0,\gamma(\pi-\psi)]\}$ with angle $\gamma(\pi-\psi)$ is covered by
 the union of the corresponding kites of the circle pattern.
\item The centers of the boundary circles lie on the boundary half lines.
\item All kites corresponding to intersecting circles are convex.
\end{enumerate}
\end{theorem}
The proof is a straight forward adaption of the proof of
Theorem~\ref{theoZgammauniq} (orthogonal case) using
Proposition~\ref{propSubharmGen2} and Lemma~\ref{lemZgammaconvex}.

\begin{theorem}
 Let $\psi\in (0,\pi)$ and $\gamma\in (0,2)\cap[\frac{\pi-2\psi}{\pi-\psi},
\frac{\pi}{\pi-\psi}]$.
Define $Z^\gamma$-circle patterns on all four sectors $\Z_\pm\times\Z_\pm$
according to Definition~\ref{defZgammapsi} and glue these patterns to a
circle pattern ${\mathscr C}_\gamma$ on a cone with cone angle $2\pi\gamma$.
Then any embedded circle pattern with the same combinatorics and intersection
angles which covers the same cone with one center of circle placed at the apex
and which has only convex kites coincides with ${\mathscr C}_\gamma$ (up to
scaling and rotation about the apex of the cone).
\end{theorem}
The proof is a simple generalization of the proof of
Theorem~\ref{theoZgammauniqpsi} to circle patterns on a cone and combinatorics
of $\Z^2$ instead of $\Z_+^2$.

The following theorem is a direct consequence of the previous theorem and
Lemmas~\ref{lem2patt} and~\ref{lemEmbeddedInfFlip} on local changes of
quasicrystallic circle patterns.

\begin{theorem}[Rigidity of quasicrystallic $Z^\gamma$-circle patterns~I]
Let $\mathscr D$ be a quasicrystallic rhombic embedding of a b-quad-graph
covering the whole plane. Let
${\cal A}=\{\pm a_1,\dots,\pm a_d\}\subset \Sp^1$ be the edge directions,
where $d>1$ and any two non-opposite elements of ${\cal A}$ are linearly
independent over $\R$. Denote by $\psi_{min}$ the minimum of the undirected
angles between any two elements of ${\cal A}$. Let $\gamma\in (0,2)$ with
$(\pi-2\psi_{min})/(\pi-\psi_{min})\leq \gamma \leq \pi/(\pi-\psi_{min})$.
Assume that the origin is a white vertex of $\mathscr D$.
Then a quasicrystallic $Z^\gamma$-circle pattern ${\mathscr
C}_\gamma$ corresponding to $\mathscr D$ and embedded on a cone with cone angle
$2\pi\gamma$ can be defined using the definition of the comparison function
$w$ on
the $2d$ sectors of $\Z^d$ which contain the quad-surface
$\Omega_{\mathscr D}$; see~\eqref{eqdefsemi} and the remarks below.
Assume further that the brick $\Pi(\Omega_{\mathscr D})$ contains the whole
lattice $\Z^d$.

Let ${\mathscr C}$ be an embedded circle pattern with the same combinatorics and
the same intersection angles which covers the same cone with one center of
circle  placed at the apex. Extend the comparison function $w$ for ${\mathscr
C}$ from $\Omega_{\mathscr D}$ to $\Z^d$.
For each $\Z^2$-sublattice which
contains two coordinate axes suppose that the corresponding circle pattern
built according to this comparison function has only convex kites. Then
${\mathscr C}$ coincides with ${\mathscr C}_\gamma$ up to
scaling and rotation about the apex of the cone.
\end{theorem}

Note that the assumption on the convexity of the kites is only a restriction
for a (small) neighborhood of the origin. This is due to Lemma~\ref{corSn}
which implies that the ratio of the radii is almost
one and thus the corresponding angles are almost the same in the isoradial
case if the combinatorial distance to the origin is big enough.

If all intersection angles of the labelling $\alpha$ are larger than $\pi/2$,
then the kites of any corresponding circle pattern are convex.
Examples of such rhombic embeddings are suitable regular hexagonal patterns.
Hexagonal circle patterns and in particular analogs of
the holomorphic mappings $z^\gamma$ have been studied by Bobenko and Hoffmann
in~\cite{BH03}.

\begin{theorem}[Rigidity of quasicrystallic $Z^\gamma$-circle patterns II]
Let $\mathscr D$ be a quasicrystallic rhombic embedding of a b-quad-graph.
Assume that the corresponding labelling $\alpha:F({\mathscr D})\to[\pi/2,\pi)$
has only values larger than $\pi/2$. Assume further that the origin is a white
vertex. Let $\gamma\in (0,2)$.

Define a quasicrystallic $Z^\gamma$-circle pattern ${\mathscr
C}_\gamma$ for $\mathscr D$ and $\alpha$ which is embedded on a cone with cone
angle $2\pi\gamma$ using the definition of the comparison function on
the $2d$ sectors of $\Z^d$ whose union contains the quad-surface
$\Omega_{\mathscr D}$; see~\eqref{eqdefsemi}. In particular, the
circle corresponding to the origin is centered at the apex of the cone.

Let ${\mathscr C}$ be an embedded circle pattern for $\mathscr D$ and
$\alpha$ which covers the same cone. Suppose that the circle corresponding to
the origin is centered
at the apex and that ${\mathscr C}$ has only convex kites. Then ${\mathscr C}$
coincides with ${\mathscr C}_\gamma$ up to scaling and rotation about the apex
of the cone.
\end{theorem}
In this case that all intersection angles are
larger than $\pi/2$ the proof of Theorem~\ref{theoZgammauniq} can be
directly adapted using Proposition~\ref{propSubharmGen2} and the following
lemma.

\begin{lemma}\label{lemsimplewalk}
Let $\mathscr D$ be a quasicrystallic rhombic embedding of a b-quad-graph which
covers the whole
plane $\C$. Let $G$ be the associated infinite graph built from white vertices.
Then the simple random walk on $G$ is recurrent.
\end{lemma}
\begin{proof}
Recall
that for the simple random walk all edges have conductance $c(e)=1$ and also
resistance $r(e)=1/c(e)=1$. Our aim is to prove that the effective resistance
$R_{\text{eff}}$ of the network $G$ with these unit conductances is infinite.
If two incident vertices of $G$ are identified, that is
the conductance of the connecting edge is increased to $\infty$ and the
resistance decreased to $0$, then the effective resistance of the new network
is certainly smaller. More generally, the proceedure of identifying a set of
vertices of $G$ is called {\em shortening} and reduces the effective
resistance; see~\cite[Section 2.2.2]{DS84} or~\cite[Theorem (2.19)]{Woe00}.
Therefore, it is sufficient to show that we have infinite effective resistance
$R_{\text{eff}}'=\infty$ for
a network $G'$ which is obtained from $G$ by shortening.

 Without loss of generality we assume that all edges of $\mathscr D$ have
length one. Since $\mathscr D$ is quasicrystallic rhombic embedding the areas
of the rhombi are uniformly bounded. Denote the lower bound by $C_1>0$.
Let $v_0\in V(G)$ be any vertex and set
$\varrho_k=4k$ for $k\in\N$. Denote by $V_k\subset V(G)$ the set of
vertices which are contained in the annulus $A_k=\{z\in\C : \varrho_{k-1}\leq
|z-v_0|< \varrho_k\}$
for $k\geq 1$. Then $V(G)=\cup_{k=1}^\infty V_k$. Identify the vertices of
each
$V_k$ to one new vertex $v_k$. Then by construction, $v_k$ is only incident to
$v_{k-1}$ for $k\geq 2$ and to $v_{k+1}$ for $k\geq 1$;
see Figure~\ref{figShortnet}.
\begin{figure}[tb]
\begin{center}
 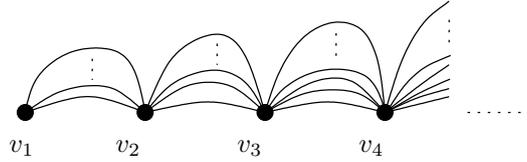
\end{center}
 \caption{An illustration of the network obtained from $G$ by
shortening.}\label{figShortnet}
\end{figure}
Denote by $|E_k|$ the number of edges which are incident to $v_k$ and
$v_{k+1}$ for $k\geq 1$. Then the effective resistance of the shortened
network
is $R_{\text{eff}}'=\sum_{k=1}^\infty 1/|E_k|$.
Furthermore $|E_k|\leq$ area of $\{z\in\C : \varrho_{k}-2\leq |z|<
\varrho_k+2\}/C_1$.
Now a simple estimation shows that $R_{\text{eff}}'=\infty$.
\end{proof}

\begin{acknowledgement}
 The author is particularly grateful to A.~I.~Bobenko for introducing her to
quasicrystallic and to $Z^\gamma$-circle patterns and for valuable discussion
and advice.
Furthermore the author thanks Veronika Schreiber for the possibility to
create figures using her computer program.

The research was partially supported by the DFG Research Unit
"Polyhedral Surfaces" and by the DFG Research Center {\sc Matheon} in
Berlin.
\end{acknowledgement}

%%%%%%%% Appendix %%%%%%%%%%%%%%%%%%
\begin{appendix}

\section{Appendix: Proof of Lemma~\ref{lemZgammaconvex}}\label{app1}

First note that kites with intersection angle larger than~$\pi/2$ are always
convex.
The remaining claims on convexity are consequences of
Proposition~3 of~\cite{Agpsi}, namely
\begin{multline}\label{eqZgammaR2psi}
\scriptstyle
(N+M)(R(z)^2-R(z+1)R(z-i) -\cos\psi R(z)(R(z-i)-R(z+1)))(R(z+i)+R(z+1)) \\
\scriptstyle
+(M-N)(R(z)^2-R(z+i)R(z+1) -\cos\psi
R(z)(R(z+i)-R(z+1))) (R(z+1)+R(z-i))=0
\end{multline}
for $z\in \V\setminus\{N+iN|N\in\N\}$, and the results of
Theorem~\ref{theoZgammapsi}.

In particular, for $\pi-\psi\leq \pi/2$ the kites with intersection angle
$\psi\leq \pi/2$ at black vertices are always convex.
This follows from equation~\eqref{eqZgammaR2psi} and
inequality~\eqref{eqRotDirectionpsi} by simple calculations.
This shows the claim for $\psi\geq \pi/2$.

If $\psi< \pi/2$ and $(\pi-2\psi)/(\pi-\psi)\leq \gamma
\leq \pi/(\pi-\psi)$,
inequality~\eqref{eqRotDirectionpsi} only implies that
for all kites with white vertices $z$ and $z+i$
and intersection angle $\psi$ the angle at the point corresponding to $z$ for
$0<\gamma<1$ and to $z+i$ for $1<\gamma<2$ respectively is smaller than $\pi$.
This excludes some types of non-convex kites, but not all.

For fixed $n>0$, let $\Gamma_n$ be the piecewise linear
curve formed by the segments $[f_{n,m},f_{n,m+1}]$ for $m\geq 0$.
By Theorem~\ref{theoZgammapsi} and its proof in~\cite{Agpsi}, these curves are
embedded without self-intersections
and the vector ${\mathbf v}_n(m) =f_{n,m+1}-f_{n,m}$ rotates clockwise for
$0<\gamma<1$ and counterclockwise for $1<\gamma <2$ along $\Gamma_n$.

Without loss of generality, we only consider the case $1<\gamma <2$ further.
For $0<\gamma<1$ the proof is very similar.
First, consider a kite on the symmetry axis, that is with white vertices
corresponding to $iK$ and
$i(K+1)$. Then the assumption $\gamma (\pi-\psi) \leq \pi$ and the properties
of the curves $\Gamma_n$ imply that
the kites on the symmetry axis are convex.
Next, we consider the intersection angles $\alpha$ and $\beta$ of the line in
direction of the vector ${\mathbf v}_n(m)$
with the oriented half lines $\R^+$ and
$\text{e}^{i\gamma(\pi-\psi)/2}\R^+$ respectively.
We additionally assume that $n$ is odd.
As the kites on the symmetry axis are convex and $\gamma (\pi-\psi) \leq \pi$,
we deduce $\alpha\leq \psi$ and $\pi/2-\psi\leq \beta\leq \pi/2$.

Finally, consider a kite with white vertex $f_{n,m}$ for odd $n$ and $m<n$. We
estimate the angles $\alpha_1$ and $\alpha_2$ at this vertex of the kites
containing the points $f_{n,m-1}$, $f_{n,m}$, $f_{n+1,m}$ and
$f_{n+1,m}$, $f_{n,m}$, $f_{n,m+1}$ respectively. Note that these kites both
have intersection angles $\psi$.
Using the above estimations
on $\alpha$ and $\beta$ and the properties of the curve $\Gamma_n$
we obtain $\pi-2\psi\leq \alpha_1$ and $\alpha_2\leq \pi$.
Therefore the angle at $f_{n-1,m-1}$ of the kite containing the points
$f_{n,m-1}$, $f_{n,m}$, $f_{n+1,m}$ is $2\pi-2\psi-\alpha_1\leq \pi$. So this
kite is convex.
Furthermore, we deduce that the kite containing the points $f_{n+1,m}$,
$f_{n,m}$, $f_{n,m+1}$ is also convex.
Consequently all kites containing the white vertex $f_{n,m}$ are convex.

Thus the $Z^\gamma$-circle pattern with $\psi< \pi/2$ and
$(\pi-2\psi)/(\pi-\psi)\leq \gamma
\leq \pi/(\pi-\psi)$ only contains convex kites.
\end{appendix}

{\small
\bibliography{dissertation,paperbibliography}
}

\end{document}

%% file: intersectionAngleKite.pdf_t
\begin{picture}(0,0)%
\includegraphics{intersectionAngleKite.pdf}%
\end{picture}%
\setlength{\unitlength}{4144sp}%
\begingroup\makeatletter\ifx\SetFigFont\undefined%
\gdef\SetFigFont#1#2#3#4#5{%
  \reset@font\fontsize{#1}{#2pt}%
  \fontfamily{#3}\fontseries{#4}\fontshape{#5}%
  \selectfont}%
\fi\endgroup%
\begin{picture}(1821,1176)(2673,-1874)
\put(3374,-1131){\makebox(0,0)[lb]{\smash{{\SetFigFont{8}{9.6}{\familydefault}{\mddefault}{\updefault}{\color[rgb]{0,0,0}$\alpha$}%
}}}}
\put(3357,-924){\makebox(0,0)[lb]{\smash{{\SetFigFont{8}{9.6}{\familydefault}{\mddefault}{\updefault}{\color[rgb]{0,0,0}$\beta$}%
}}}}
\end{picture}%

%% file: ExampleShort2.pdf_t
\begin{picture}(0,0)%
\includegraphics{ExampleShort2.pdf}%
\end{picture}%
\setlength{\unitlength}{4144sp}%
\begingroup\makeatletter\ifx\SetFigFont\undefined%
\gdef\SetFigFont#1#2#3#4#5{%
  \reset@font\fontsize{#1}{#2pt}%
  \fontfamily{#3}\fontseries{#4}\fontshape{#5}%
  \selectfont}%
\fi\endgroup%
\begin{picture}(3105,972)(1156,-3271)
\put(1171,-3211){\makebox(0,0)[lb]{\smash{{\SetFigFont{10}{12.0}{\familydefault}{\mddefault}{\updefault}{\color[rgb]{0,0,0}$v_1$}%
}}}}
\put(1801,-3211){\makebox(0,0)[lb]{\smash{{\SetFigFont{10}{12.0}{\familydefault}{\mddefault}{\updefault}{\color[rgb]{0,0,0}$v_2$}%
}}}}
\put(2521,-3211){\makebox(0,0)[lb]{\smash{{\SetFigFont{10}{12.0}{\familydefault}{\mddefault}{\updefault}{\color[rgb]{0,0,0}$v_3$}%
}}}}
\put(3241,-3211){\makebox(0,0)[lb]{\smash{{\SetFigFont{10}{12.0}{\familydefault}{\mddefault}{\updefault}{\color[rgb]{0,0,0}$v_4$}%
}}}}
\end{picture}%